\documentclass{article}
\usepackage[utf8]{inputenc}
\usepackage[T1]{fontenc}
\usepackage[english]{babel}

\usepackage{mathtools}
\usepackage{bm}
\usepackage{dsfont}
\usepackage{booktabs}

\usepackage{graphicx}
\usepackage{float}
\usepackage{natbib}
\usepackage{amsmath,amssymb,amsthm,mathtools}

\usepackage{multirow}
\usepackage{xcolor}

\usepackage{enumitem}

\usepackage[doublespacing]{setspace}
\usepackage{url}
\usepackage[
  colorlinks=true,
  citecolor=black,
  linkcolor=blue,
  urlcolor=blue,
  bookmarksnumbered=true
]{hyperref}
\usepackage{cleveref}
\newtheorem{theorem}{Theorem}
\numberwithin{theorem}{section}

\numberwithin{proposition}{section}

\newtheorem{lemma}{Lemma}
\numberwithin{lemma}{section}

\newtheorem{remark}{Remark}

\newtheorem{hyp}{H}

\crefname{theorem}{theorem}{theorems}
\crefname{lemma}{lemma}{lemmas}
\crefname{proposition}{proposition}{propositions}
\crefname{definition}{definition}{definitions}
\crefname{example}{example}{examples}
\crefname{remark}{remark}{remarks}
\crefname{hyp}{H}{H}

\Crefname{theorem}{Theorem}{Theorems}
\Crefname{lemma}{Lemma}{Lemmas}
\Crefname{proposition}{Proposition}{Propositions}
\Crefname{definition}{Definition}{Definitions}
\Crefname{example}{Example}{Examples}
\Crefname{remark}{Remark}{Remarks}
\Crefname{hyp}{H}{H}

\title{Nonparametric estimation of a common marginal hazard function under censoring and spatio-temporal dependence}

\author{
Ibrahim Loudy\textsuperscript{a}, 
		Kowir Pambo Bello\textsuperscript{b}, 
		St\'ephane Bouka\textsuperscript{a}\thanks{Corresponding author: St\'ephane Bouka. Email: stephane.bouka@univ-mauku.org}\, and\\
		Guy Martial Nkiet\textsuperscript{a}}
	\date{\textsuperscript{a}Laboratoire de Probabilit\'es, Statistique et Informatique,
	Unit\'e de Recherche de Math\'ematiques et Informatique,
	Universit\'e des Sciences et Techniques de Masuku, Franceville, Gabon\\
	\textsuperscript{b}Laboratoire de Recherche en Math\'ematiques,
	\'Ecole Normale Sup\'erieure, Libreville, Gabon}

\begin{document}


\maketitle

\begin{abstract}
	This article investigates the nonparametric estimation of a common pooled marginal hazard function from right-censored durations exhibiting spatio-temporal dependence.
	We propose two estimators based respectively on the Nadaraya--Watson method and the \(k\)-nearest neighbors (\(k\)-NN) approach. Under assumptions of stationarity, regularity, non-informative censoring, and weak spatio-temporal dependence, we establish their pointwise and uniform convergence rates as well as their optimal asymptotic order. A simulation study compares their finite-sample performance.
	Two real-data applications with right censoring are considered: a multicenter clinical trial with 53.8\,\% censoring and 1,746 geolocated and time-stamped requests submitted to New York City's 311 municipal service, of which 35.5\,\% are censored.
	The empirical results illustrate the relevance of the proposed methods in the presence of spatial and temporal dependence.
\end{abstract}



	\section{Introduction}\label{sec1}
	\noindent
	Statistical problems related to the modeling and analysis of spatio-temporal data have attracted increasing attention in recent years. Understanding phenomena that evolve simultaneously over time and space has become essential in many fields of application, including meteorology, epidemiology, ecology, geography, economics, and urban planning. Such phenomena often exhibit complex dependence structures that cannot be adequately captured by purely temporal or purely spatial models.
	In nonparametric statistics, however, the modeling of such data remains relatively underexplored. Recent contributions by \cite{zhou2025machine}, \cite{zhang2025multidimensional}, and \cite{loudy2025non} illustrate this trend by proposing nonparametric approaches adapted to the spatio-temporal framework.
	The estimation of the hazard function plays a central role in survival analysis. Although this problem has been widely studied for independent or temporal observations (see, in particular, \cite{jiang11}, \cite{mohammed19}, \cite{douma18}, and \cite{mokkadem21}), the case of censored data exhibiting spatio-temporal dependence remains underdeveloped.
	Nevertheless, some contributions in the spatial and spatio-temporal literature are worth mentioning. For instance, \cite{li07} proposes a nonparametric estimator of the hazard function for a spatial random field and establishes its asymptotic normality in order to construct confidence intervals. Similarly, \cite{laksaci10} develops a kernel estimator of the conditional hazard function in a functional space and proves its almost sure convergence for spatially dependent data.
	Some works address the issue more indirectly through the study of relative risk or occurrence intensity. For example, \cite{park2023exploring} analyze the spatial distribution of restaurant closures in Seoul during the COVID-19 pandemic using a relative risk surface; \cite{zhang11} propose stKDE, a method that jointly estimates disease risk in space and time; and \cite{hazelton2023shrinkage} introduces a Lasso-type estimator of spatial relative risk between cases and controls to better identify high-risk areas despite data inhomogeneity.
	Our objective is to estimate a common pooled marginal hazard function $\lambda(x)$, identical for all sites $\mathbf{s}$ and all periods $t$, under spatio-temporal dependence among the observations.
	We propose two nonparametric estimators of this common marginal hazard function, constructed respectively from the Nadaraya--Watson method and the \(k\)-nearest neighbors (\(k\)-NN) approach. For both estimators, we establish pointwise and uniform convergence rates, and then compare their behavior through a simulation study and two real-data applications.
	The article is organized as follows. Section~\ref{sec2} presents the proposed estimators. Section~\ref{sec3} introduces the assumptions as well as the main theoretical results. Section~\ref{sec4} is devoted to the simulation study and the real-data applications. Section~\ref{sec5} concludes the paper, while Section~\ref{sec6} gathers the proofs.

	\section{\texorpdfstring{Estimation of a common marginal hazard function under censoring and spatio-temporal dependence}{Estimation of a common marginal hazard function under censoring and spatio-temporal dependence}}
	\label{sec2}
	
	Throughout, we write $\mathbb{N}^* = \{1,2,3,\ldots\}$ for the set of strictly positive integers and $\mathbb{R}_+^* = (0,+\infty)$ for the set of strictly positive reals.
	We consider a spatio-temporal survival process observed over $T \in \mathbb{N}^*$ temporal periods indexed by
	$
	t=1,\dots,T,
	$
	and on a $d$-dimensional spatial grid ($d \geqslant 1$)
	$
	\Lambda_{\mathbf n},
	$
	where $\mathbf{n} = (n_1,\ldots,n_d) \in (\mathbb{N}^*)^d$ collects the number of observation sites along each spatial dimension.
	For each time $t$ and each spatial site $\mathbf i\in\Lambda_{\mathbf n}$, we denote by
	$
	\mathcal{T}_{t,\mathbf{i}}\in\mathbb R_{+}
	$
	the latent survival time of interest, that is, the duration between the origin of observation and the occurrence of the event under study.
	When the data are right-censored, we also introduce a censoring time
	$
	C_{t,\mathbf i},
	$
	and we only observe
	\[
	X_{t,\mathbf{i}}=\min(\mathcal{T}_{t,\mathbf{i}},C_{t,\mathbf i}),
	\qquad
	\delta_{t,\mathbf i}=\mathbf 1_{\{\mathcal{T}_{t,\mathbf{i}}\leqslant C_{t,\mathbf i}\}},
	\]
	where $\delta_{t,\mathbf i}$ indicates whether the observed duration is uncensored.
	The available data are thus given by the family
	\[
	\bigl\{(X_{t,\mathbf{i}},\delta_{t,\mathbf i}),\; 1\leqslant t\leqslant T,\; \mathbf i\in\Lambda_{\mathbf n}\bigr\}.
	\]
	In this framework, it is important to clearly distinguish between two notions of time.
	On the one hand, the index $t$ only identifies the observation period and does not represent a survival duration.
	On the other hand, the variable $\mathcal{T}_{t,\mathbf{i}}$ is a positive real duration measuring the elapsed time until the occurrence of the event of interest.
	For example, if $t=3$, this simply means that the observation is collected during the third study period, whereas the equality
	$
	\mathcal{T}_{3,\mathbf{i}}=12.4
	$
	indicates that, at site $\mathbf i$ during this period, the event occurs after $12.4$ time units.
	Such a structure arises in many applications.
	For instance, if $\mathbf i$ denotes a hospital and $t$ a year of observation, then $\mathcal{T}_{t,\mathbf{i}}$ may represent the time to relapse of a patient followed in this hospital during year $t$.
	Similarly, if $\mathbf i$ denotes an industrial site and $t$ a quarter, then $\mathcal{T}_{t,\mathbf{i}}$ may represent the time to failure of a device observed at this site during the corresponding quarter.
	Let $F$ and $f$ denote the distribution function and density of the survival variable $\mathcal{T}_{t,\mathbf{i}}$, $G$ and $g$ those of the censoring time $C_{t,\mathbf i}$, and $H$ and $\varphi_X$ the distribution function and density of the observed variable
	$
	X_{t,\mathbf{i}}=\min(\mathcal{T}_{t,\mathbf{i}},C_{t,\mathbf i}).
	$
We assume that the processes
$
(\mathcal{T}_{t,\mathbf{i}})_{\mathbf i\in\Lambda_{\mathbf n},\, 1\le t\le T}
$
and
$
(C_{t,\mathbf i})_{\mathbf i\in\Lambda_{\mathbf n},\, 1\le t\le T}
$
are $\tau$-weakly dependent in the sense of \cite{dedecker04}, and that they are independent of each other.
The common marginal hazard function associated with $\mathcal{T}_{t,\mathbf{i}}$ is defined by
	\[
	\lambda(x)=\frac{f(x)}{1-F(x)}.
	\]
	Under right censoring, we have
	\[
	1-H(x)=\mathbb P(X>x)=\mathbb P(\mathcal{T}>x,\,C>x).
	\]
	By independence between $\mathcal{T}$ and $C$, it follows that
	\[
	1-H(x)=(1-F(x))(1-G(x)),
	\]
	which yields
	\[
	1-F(x)=\frac{1-H(x)}{1-G(x)},
	\qquad \text{for } 1-G(x)>0.
	\]
	Substituting this expression into the definition of the hazard function, we obtain
	\[
	\lambda(x)=\frac{f(x)}{1-F(x)}
	=\frac{f(x)\bigl(1-G(x)\bigr)}{1-H(x)}.
	\]
	By setting
	\[
	f^{*}(x)=f(x)\bigl(1-G(x)\bigr),
	\]
	we can write
	\[
	\lambda(x)=\frac{f^{*}(x)}{1-H(x)},
	\qquad x\in\mathcal Z,
	\]
	where
	\[
	\mathcal Z
	=\left]-\infty,\overline{x}\right],
	\qquad
	\overline{x}=\sup\{x\in\mathbb R:\,1-H(x)>0\}.
	\]
	This representation is simply a rewriting of the classical hazard function $\lambda(x) = f(x)/(1-F(x))$ in terms of the observable quantities $f^*$ and $H$, which are estimable from the censored data. Indeed, when $G(x)=0$ (no censoring), we have $f^*(x)=f(x)$ and $1-H(x)=1-F(x)$, so that the two expressions coincide:
	$
	\lambda(x)=f(x)/(1-F(x)).
	$
	To estimate this common marginal hazard function from censored observations exhibiting spatio-temporal dependence, we propose two kernel-based estimators.
	The first relies on a deterministic smoothing window of Nadaraya--Watson type, while the second is based on a random window defined through the $k$-nearest neighbors method.
	The Nadaraya--Watson-type estimator is defined by
	\[
	\widetilde{\lambda}_{T,\mathbf n}(x)
	=
	\frac{\widetilde{f}^{*}_{T,\mathbf n}(x)}
	{1-\widetilde{H}_{T,\mathbf n}(x)},
	\]
	where
	\[
	\widetilde{f}^{*}_{T,\mathbf n}(x)
	=
	\frac{1}{\hat{\mathbf{n}}\,T\,\widetilde h_{T,\mathbf n}}
	\sum_{\mathbf i\in\Lambda_{\mathbf n}}\sum_{t=1}^{T}
	\delta_{t,\mathbf i}\,
	K\!\left(\frac{x-X_{t,\mathbf{i}}}{\widetilde h_{T,\mathbf n}}\right),
	\]
	and
	\[
	\widetilde{H}_{T,\mathbf n}(x)
	=
	\frac{1}{\hat{\mathbf{n}}\,T}
	\sum_{\mathbf i\in\Lambda_{\mathbf n}}\sum_{t=1}^{T}
	L\!\left(\frac{x-X_{t,\mathbf{i}}}{\widetilde h_{T,\mathbf n}}\right).
	\]
	Here,
	$
	\hat{\mathbf{n}}=n_1\times\cdots\times n_d
	$ denotes the total number of spatial locations in the observation grid.
	$K:\mathbb R\to\mathbb R$ is a continuous, positive, and integrable kernel, and
	\[
	L(u)=\int_{-\infty}^{u}K(v)\,dv.
	\]
	The quantity $\widetilde h_{T,\mathbf n}$ denotes a deterministic smoothing parameter chosen to balance bias and variance.
	The second estimator is a kernel estimator based on the $k$-nearest neighbors method, defined by
	\[
	\widehat{\lambda}_{T,\mathbf n}(x)
	=
	\frac{\widehat{f}^{*}_{T,\mathbf n}(x)}
	{1-\widehat{H}_{T,\mathbf n}(x)},
	\]
	with
	\[
	\widehat{f}^{*}_{T,\mathbf n}(x)
	=
	\frac{1}{\hat{\mathbf{n}}\,T\,h_{T,\mathbf n}(x)}
	\sum_{\mathbf i\in\Lambda_{\mathbf n}}\sum_{t=1}^{T}
	\delta_{t,\mathbf i}\,
	K\!\left(\frac{x-X_{t,\mathbf{i}}}{h_{T,\mathbf n}(x)}\right),
	\]
	and
	\[
	\widehat{H}_{T,\mathbf n}(x)
	=
	\frac{1}{\hat{\mathbf{n}}\,T}
	\sum_{\mathbf i\in\Lambda_{\mathbf n}}\sum_{t=1}^{T}
	L\!\left(\frac{x-X_{t,\mathbf{i}}}{h_{T,\mathbf n}(x)}\right).
	\]
	In this case, the smoothing window $h_{T,\mathbf n}(x)$ is random and depends on the target point $x$ and the number of selected neighbors.
	It is defined by
	\begin{equation}\label{eq:knn_bandwidth}
		h_{T,\mathbf n}(x)
		=
		\min\left\{
		h\in\mathbb R_{+}^{*}:
		\sum_{\mathbf i\in\Lambda_{\mathbf n}}\sum_{t=1}^{T}
		\mathbf 1_{\{|X_{t,\mathbf{i}}-x|\leqslant h\}}
		=
		k_{T,\mathbf n}
		\right\},
	\end{equation}
	where $(k_{T,\mathbf n})$ is an increasing sequence of strictly positive integers.
	Both estimators rely on the same decomposition of the hazard function into a numerator estimating $f^{*}$ and a denominator estimating the survival function of the observed variable.
	Their main difference lies in the choice of the smoothing window: it is fixed in the Nadaraya--Watson approach, whereas it adapts locally to the density of observations in the $k$-nearest neighbors approach.
The window \(h_{T,\mathbf n}(x)\), defined by the distance to the \(k_{T,\mathbf n}\)-th nearest neighbor, adapts to the local density of the observations \cite{ahmed23,mintogo2025rates}, thereby providing adaptive smoothing.
	These two strategies therefore provide complementary tools for the nonparametric estimation of a common marginal hazard function from right-censored data exhibiting spatio-temporal dependence.

	\section{Assumptions and Results}
	\label{sec3}
	In this section, we first define the notion of weak dependence that will be used thereafter, then we specify the assumptions involved in the main results before stating them.
	
	\subsection{Weak dependence}
	
	In this paper, we use the notion of weak dependence formulated by \cite{amiri18} and \cite{loudy2025non}.
	More precisely, given a function \(\psi : \mathbb{N}^2 \times \mathbb{R}_+^2 \to \mathbb{R}_+\), a spatio-temporal process $\{X_{t,\mathbf{s}}, \mathbf{s} \in \mathcal{S}, t \in \{1,\ldots,T\}\}$ is said to be \(\psi\)-weakly dependent if there exists a sequence \(\{\nu(r,m), (r,m) \in \mathbb{N}^2\}\), decreasing (with respect to \(r\) and \(m\)) to $0$ as \(r, m \to \infty\), such that, for all subsets of locations
	$\mathcal{S}_1 = \{\mathbf{s}'_1,\cdots,\mathbf{s}'_u\}$ and $\mathcal{S}'_1 = \{\mathbf{s}_1,\cdots,\mathbf{s}_v\}$ of $\mathcal{S}$, for any collection of $u$ and $v$ time instants
	$(t'_1,\cdots,t'_u)$ and $(t_1,\cdots,t_v)$ with $t'_1 \leqslant t'_2 \leqslant \cdots \leqslant t'_u  + z = t_1 \leqslant \cdots \leqslant t_v$, and for all bounded functions $g_0$ and $g_1$, the following inequality holds:
	\begin{align*}
		&|cov(g_0(X_{t'_1,\textbf{s}'_{1}},\cdots, X_{t'_u,\textbf{s}'_{u}}), g_1(X_{t_1,\textbf{s}_1},\cdots, X_{t_v,\textbf{s}_{v}}))|\\
		&\hspace{1cm}\leqslant \psi\left(u\Vert g_1\Vert_\infty, v\Vert g_0\Vert_\infty, Lip(g_0), Lip(g_1)\right) \nu(z,w),
	\end{align*}
	where $w := dist(\mathcal{S}_1, \mathcal{S}'_1)$ is the distance between the subsets $\mathcal{S}_1$ and $\mathcal{S}'_1$, and
	\[
	\operatorname{Lip}(g) := \sup_{x \ne y}\frac{|g(x)-g(y)|}{\Vert x-y\Vert_{1}} < \infty,
	\]
	with $\Vert u\Vert_1=\sum_i\vert u_i\vert$.
	When $\psi(j, k, u, v) = ju + kv$, the coefficient $\nu$ is denoted by $\eta$, and in this case the process is simply said to be $\eta$-weakly dependent \cite{kallabis06}.
	As in \cite{amiri18}, in order to obtain uniform consistency results, we will use a spatio-temporal counterpart of the notion of $\tau$-weak dependence introduced by \cite{dedecker04}, the $\eta$-weak dependence being more general but insufficient to obtain uniform consistency results.
	Let $\mathcal{L}(\mathbb{R})$ be the set of functions $f: \mathbb{R} \to \mathbb{R}$ such that $\mathrm{Lip}(f) < 1$.
	For any integrable real-valued random variable $X$, defined on a probability space $(\Omega, \mathcal{A}, \mathbb{P})$, and for any $\sigma$-algebra $\mathcal{M}$ included in $\mathcal{A}$, we set
	\[
	\tau(\mathcal{M}, X) = \mathbb{E} \left[ \sup_{f \in \mathcal{L}(\mathbb{R})} \left| \mathbb{P}_{X|\mathcal{M}}(f) - \mathbb{P}_X(f) \right| \right],
	\]
	where $\mathbb{P}_X$ and $\mathbb{P}_{X|\mathcal{M}}$ denote respectively the distribution of $X$ and the conditional distribution of \(X\) given \(\mathcal{M}\).
	We define
\[
\begin{aligned}
	\tau(z,w)
	&= {}\\
	&\hspace{-0.7cm}\sup_{n>0}\max_{1\leqslant v\leqslant n}\frac{1}{v}
	\sup_{\substack{
			\mathcal{S}_1=\{s'_1,\ldots,s'_u\}\subset \mathcal{S},\\
			\mathcal{S}'_1=\{s_1,\ldots,s_v\}\subset \mathcal{S},\\
			w=\operatorname{dist}(\mathcal{S}_1,\mathcal{S}'_1),\\
			t'_1\leqslant t'_2\leqslant \cdots \leqslant t'_u+z=t_1\leqslant \cdots \leqslant t_v
	}}
	\tau\!\left[
	\sigma\!\left(X_{t'_1,s'_1},\ldots,X_{t'_u,s'_u}\right),
	X_{t_1,s_1},\ldots,X_{t_v,s_v}
	\right],
\end{aligned}
\]
	where $\operatorname{dist}(\mathcal{S}_1, \mathcal{S}'_1)$ is the distance between the subsets $\mathcal{S}_1$ and $\mathcal{S}'_1$.
	We shall say that $\{X_{t,\mathbf{s}}, \mathbf{s} \in \mathcal{S}, t \in \{1,\ldots,T\}\}$ is $\tau$-weakly dependent if $\tau(z, w) \to 0$ as $z, w \to \infty$.
	
	\subsection{Assumptions}
	
	\noindent From now on, we shall say that \((T, \widehat{\textbf{n}}) \to \infty\) whenever
	\[
	\min_{1  \leqslant  \ell  \leqslant  d} \{n_\ell\} \to \infty, \quad T \to \infty, \quad \text{and} \quad \frac{n_j}{n_\ell}  \leqslant  c_1 < \infty \quad \text{for all } 1  \leqslant  j, \ell  \leqslant  d.
	\]
	\noindent We make the following assumptions:
	
\begin{hyp}\label{H1}
Let $\mathbb{T}=\{1,\ldots,T\}$. The spatio-temporal process
$
\{X_{t,\mathbf{i}},\, \mathbf{i}\in\Lambda_{\mathbf{n}},\, t\in\mathbb{T}\}
$
is $\tau$-weakly dependent. Moreover, there exist two functions $\tau_1$ and $\tau_2$ satisfying
$
\tau_1(0)=\tau_2(0)=1,
$
such that, for all $z,w\in\mathbb{R}$,
$
\tau(z,w)=\tau_1(z)\tau_2(w),
$
with
$
\tau_1(z)=O\bigl(\exp(-\rho_1 z)\bigr)
\quad \text{and} \quad
\tau_2(w)=O\bigl(\exp(-\rho_2 w)\bigr),
$
where $\rho_1>3/2$ and $\rho_2>3/2$.
\end{hyp}

\begin{hyp}\label{H2}
The kernel \(K\) is a bounded, Lipschitz, nonnegative density with compact support \(D=[-a,a]\), with $a>0$, and symmetric. Moreover, for every \(u\in\mathbb R\) and every \(t\in]0,1[\), \(K(u)\leqslant K(tu)\).
\end{hyp}

\begin{hyp}\label{H3}
$f^*$ is bounded, belongs to $C^2(\mathbb{R})$ with bounded derivatives; moreover, the density $\varphi_X$ of the observed variable $X$ belongs to $C^1(\mathbb{R})$ and has a bounded derivative.
\end{hyp}

\begin{hyp}\label{H4}
$k_{T, \mathbf{n}} \to \infty$; $\dfrac{k_{T, \mathbf{n}}}{\widehat{\mathbf{n}} T} \to 0$ and $\dfrac{k_{T, \mathbf{n}}}{\log(\widehat{\mathbf{n}} T)} \to \infty$ as $(T,\widehat{\textbf{n}}) \to +\infty$.
\end{hyp}

\begin{hyp}\label{H5}
$\widetilde{h}_{T, \mathbf{n}} \to 0$ and $\widehat{\textbf{n}}T\widetilde{h}_{T, \mathbf{n}} \to \infty$ as $(T,\widehat{\textbf{n}}) \to \infty$.
\end{hyp}

\begin{hyp}\label{H6}
The process $\{X_{t, \mathbf{i}}\}$ is strictly stationary, and the joint density $\varphi_{(X_{t, \mathbf{i}},\, X_{s, \mathbf{j}})}$ of the pair $(X_{t, \mathbf{i}},\, X_{s, \mathbf{j}})$ is such that
$
M := \displaystyle\sup_{(t,\textbf{i}) \ne (s ,\textbf{j})}\Vert \varphi_{(X_{t,\textbf{i}} , X_{s , \textbf{j}})} - \varphi_X \otimes \varphi_X\Vert_{\infty} < \infty.
$
\end{hyp}

\begin{hyp}\label{H7}
For every compact set $\mathcal C$ included in $\mathcal Z$, there exists a constant $c_{\mathcal C}>0$ such that
	$
	\inf_{x\in\mathcal C}\{1-H(x)\}\geqslant c_{\mathcal C},
	$
\end{hyp}

\begin{remark}
	\begin{itemize}[leftmargin=2em,label=\textbullet]
		\item {\hyperref[H1]{H1}} describes the dependence structure. The qualitative condition $\tau(z,w)\to0$ defines weak dependence \cite{dedecker04,amiri18}, a framework that may include non-mixing processes \cite{andrews84,doukhan99}. The exponential decay imposed in {\hyperref[H1]{H1}} is a sufficient condition for the summability of the coupling errors.

	\item \hyperref[H2]{H2} collects the usual kernel conditions found, in particular, in \cite{ahmed23,mintogo2025rates}. They are satisfied by the Epanechnikov, biweight, and triweight kernels and make it possible to control the random window of the $k$-NN estimator.
	
\item \hyperref[H3]{H3} imposes regularity conditions on $f^*$ and $\varphi_X$, which are used to control the bias and establish the convergence rates \cite{khardani12,lecoutre95}.

	\item \hyperref[H4]{H4} specifies the regime of $k_{T,\mathbf n}$ in order to ensure increasing local information, a shrinking window, and control of the stochastic term \cite{ahmed23,mintogo2025rates}.
	\item \hyperref[H5]{H5} specifies the regime of $\widetilde h_{T,\mathbf n}$ in order to ensure that the bias vanishes and that there is sufficient local information for consistency of the Nadaraya--Watson estimator \cite{lecoutre95,amiri18}.
	\item \hyperref[H6]{H6} imposes strict stationarity and a uniform bound on the difference between the joint density and the product of the marginal densities, a standard condition for the study of kernel estimators under spatial dependence \cite{el2020nadaraya,el2021class}.

		\item \hyperref[H7]{H7} ensures that $1-H$ remains uniformly positive on the compact sets under consideration, in order to control the quotient defining the hazard estimator and its uniform convergence \cite{liebscher2002kernel}.
	\end{itemize}
\end{remark}

	\subsection{Results}
	
\noindent The following theorem establishes the pointwise convergence rate of the Nadaraya--Watson-type estimator of the common marginal hazard function, in the context of censored weakly dependent data:
	\begin{theorem}
		Under assumptions \hyperref[H1]{H1}, \hyperref[H2]{H2}, \hyperref[H3]{H3}, \hyperref[H5]{H5}, \hyperref[H6]{H6}, for every \( x \in \mathcal{Z} \), one has:
		\[
		\left|  \widetilde{\lambda}_{T,\mathbf{n}}(x) - \lambda(x) \right| = O\left( \widetilde{h}_{T,\mathbf{n}}^{\,2} + \sqrt{\dfrac{\log( \widehat{\mathbf{n}} T)}{\widehat{\mathbf{n}} T \widetilde{h}_{T,\mathbf{n}}}} \right) \hspace{1cm} \text{a.s.}
		\]
		\label{th1}
	\end{theorem}
	
\noindent The following theorem specifies the strong uniform consistency rate of the Nadaraya--Watson-type estimator of the common marginal hazard function in the framework of censored weakly dependent data:
	\begin{theorem}
	Under assumptions \hyperref[H1]{H1}, \hyperref[H2]{H2}, \hyperref[H3]{H3}, \hyperref[H5]{H5}, \hyperref[H6]{H6}, and \hyperref[H7]{H7}, for every compact set \( \mathcal{C} \) included in \( \mathcal{Z} \),
		one has:
		\[
		\sup_{x \in \mathcal{C}} \left|  \widetilde{\lambda}_{T,\mathbf{n}}(x) - \lambda(x) \right| =  O\left( \widetilde{h}_{T,\mathbf{n}}^{\,2} + \sqrt{\dfrac{\log( \widehat{\mathbf{n}} T)}{\widehat{\mathbf{n}} T \widetilde{h}_{T,\mathbf{n}}}} \right) \hspace{1cm} \text{a.s.}
		\]
		\label{th2}
	\end{theorem}
	
\noindent The following theorem establishes the pointwise convergence rate of the \(k\)-nearest neighbors estimator of the common marginal hazard function, in the context of censored weakly dependent data.
	\begin{theorem}
	Under assumptions \hyperref[H1]{H1}, \hyperref[H2]{H2}, \hyperref[H3]{H3}, \hyperref[H4]{H4}, \hyperref[H6]{H6}, for every \( x \in \mathcal{Z} \), one has:
		\[
		\left|  \widehat{\lambda}_{T,\mathbf{n}}(x) - \lambda(x) \right| =O\left( \left( \dfrac{k_{T,\mathbf{n}}}{\widehat{\mathbf{n}} T} \right)^2 + \sqrt{ \dfrac{ \log( \widehat{\mathbf{n}} T ) }{ k_{T,\mathbf{n}} } } \right). \hspace{1cm} \text{a.s.}
		\]
		\label{th3}
	\end{theorem}
	
\noindent The following theorem gives the strong uniform consistency rate of the \(k\)-nearest neighbors estimator of the common marginal hazard function, in the context of censored weakly dependent data.
	\begin{theorem}
		Under assumptions \hyperref[H1]{H1}, \hyperref[H2]{H2}, \hyperref[H3]{H3}, \hyperref[H4]{H4}, \hyperref[H6]{H6}, and \hyperref[H7]{H7}, for every compact set \( \mathcal{C} \) included in \( \mathcal{Z} \),
		one has:
		\[
		\sup_{x \in \mathcal{C}} \left|  \widehat{\lambda}_{T,\mathbf{n}}(x) - \lambda(x) \right| = O\left( \left( \dfrac{k_{T,\mathbf{n}}}{\widehat{\mathbf{n}} T} \right)^2 + \sqrt{ \dfrac{ \log( \widehat{\mathbf{n}} T ) }{ k_{T,\mathbf{n}} } } \right). \hspace{1cm} \text{a.s.}
		\]
		\label{th4}
	\end{theorem}
	
	\begin{remark}
		\noindent
		The previous theorems show that the Nadaraya--Watson-type and \(k\)-nearest neighbors estimators achieve the same order of convergence for estimating \(\lambda(x)\). Balancing bias and variance yields the optimal choices \(\widetilde{h}_{T,\mathbf{n}} \approx \left( \frac{\log(\widehat{\mathbf{n}}T)}{\widehat{\mathbf{n}}T} \right)^{1/5}\) and \(k_{T,\mathbf{n}} \approx \left( (\widehat{\mathbf{n}}T)^4 \log(\widehat{\mathbf{n}}T) \right)^{1/5}\), leading in both cases to the rate
		\[
		O\left( \left( \frac{\log(\widehat{\mathbf{n}}T)}{\widehat{\mathbf{n}}T} \right)^{2/5} \right).
		\]
		Hence, both methods are asymptotically equivalent from a theoretical point of view, and their comparison mainly depends on practical considerations, which will be assessed through simulations.
	\end{remark}

	\begin{remark}
		The assumptions \hyperref[H1]{H1}--\hyperref[H7]{H7} are primarily formulated to obtain the stated convergence rates and uniform results; they can be weakened when only pointwise consistency is sought. In particular, in \hyperref[H1]{H1}, polynomial decay is sufficient. For \hyperref[H2]{H2}, it is enough for the kernel $K$ to be integrable, square-integrable, and to satisfy
		$
		\int_{\mathbb R}K(u)\,du=1,
		$
	 in \hyperref[H3]{H3}, continuity of the functions concerned is sufficient for consistency. Assumptions \hyperref[H4]{H4} and \hyperref[H5]{H5} can be replaced by the usual conditions
		$
		k_{T,\mathbf n}\to\infty,\quad
		\frac{k_{T,\mathbf n}}{\widehat{\mathbf n}T}\to0,
		$
		and, for Nadaraya--Watson, by $\widetilde h_{T,\mathbf n}\to0$. Finally, for \hyperref[H6]{H6} and \hyperref[H7]{H7}, it is sufficient that the covariance terms be asymptotically negligible and that $1-H(x)>0$. The strengthened versions of these assumptions are therefore mainly needed to explicitly control the bias and variance, establish uniform convergence, and obtain the convergence rates in the theorems.
	\end{remark}

	\section{\texorpdfstring{Simulations and Applications to Real Data}{Simulations and Applications to Real Data}}\label{sec4}

\subsection{\texorpdfstring{Simulation}{Simulation}}

\begin{figure}[H]
	\centering
	\includegraphics[width=0.75\textwidth]{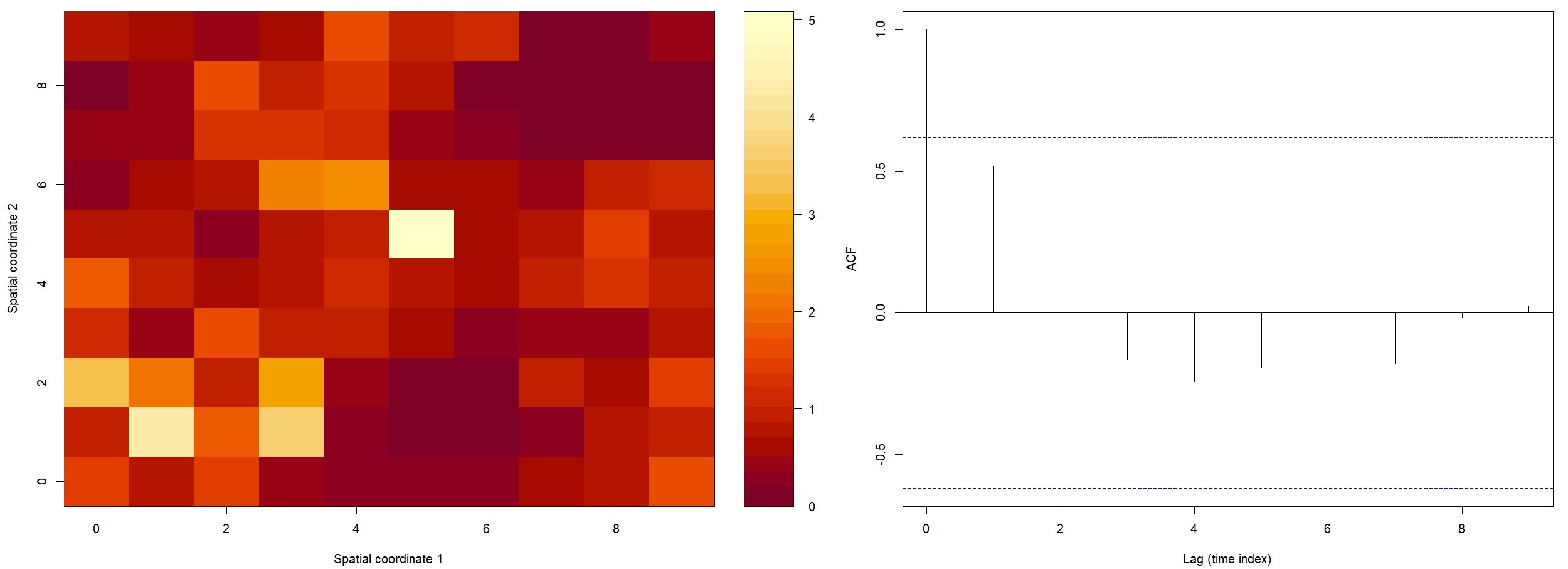}
	\caption{Simulated field under strong dependence: observations at the \(10\times10\) sites (left) and temporal autocorrelation function at site \((5,5)\) (right).}
	\label{fig:sim_domain}
\end{figure}

The simulation design is consistent with the asymptotic regime of Section~\ref{sec3}. For a configuration \(m^3\), \(m\times m\) sites are observed at \(m\) time points, giving \(N=m^3\) observations. Thus, for \(m=5,8,10\), the spatial and temporal dimensions increase jointly, corresponding to an increasing-domain asymptotic regime rather than an \emph{infill} regime in which a fixed domain is progressively densified. This choice notably avoids a situation in which, under strong spatial dependence, adding very close sites provides only a small gain in effective information. It is therefore directly consistent with the balanced asymptotic regime \((T,\mathbf n)\to\infty\) considered in the theorems of Section~\ref{sec3}.
The latent Gaussian fields associated with the event and censoring are centered, independent, and characterized by the covariance
\[
	\operatorname{Cov}\{Z_{t,\mathbf i},Z_{t',\mathbf j}\}
	=
	\exp\{-a\|\mathbf i-s\mathbf j\|\}\exp\{-b|t-t'|\}.
\]
We set \(a=b=2\) for weak dependence and
	\(a=b=0.70\) for strong dependence.
The fields are then transformed as
$
U^{\mathcal T}_{t,\mathbf i}
=
\Phi\!\left(Z^{\mathcal T}_{t,\mathbf i}\right),
\quad
U^{C}_{t,\mathbf i}
=
\Phi\!\left(Z^{C}_{t,\mathbf i}\right),
$
which produces two uniform fields that preserve the dependence induced by the Gaussian copula \cite{joe2014dependence,nelsen2006introduction}. Their independence ensures non-informative censoring \cite{kalbfleisch2002statistical,klein2003survival}.
To evaluate the estimators under contrasting hazard shapes, we consider two survival models, denoted by \(M_1\) and \(M_2\). The survival times are generated by inversion,
$
\mathcal T_{t,\mathbf i}
=
F_{M_w}^{-1}\!\left(U^{\mathcal T}_{t,\mathbf i}\right),
\quad w\in\{1,2\},
$
where \(F_{M_w}\) denotes the parametric distribution function of model \(M_w\): \(M_1\) is a Weibull distribution with parameters \((0.7,1)\), and \(M_2\) is a log-logistic distribution with parameters \((1,2)\), with
\[
\lambda_1(x)=0.7x^{-0.3},
\qquad
\lambda_2(x)=\frac{2x}{1+x^2}.
\]
The censoring times are generated from an exponential distribution,
$
C_{t,\mathbf i}
=
F_{\mathrm{Exp}}^{-1}\!\left(U^C_{t,\mathbf i};\mu\right),
$
where \(\mu\) is calibrated to obtain \(p_c\in\{20\%,60\%\}\). The selected values are \(0.224\) and \(1.869\) for \(M_1\), and \(0.170\) and \(0.926\) for \(M_2\). Finally,
$
X_{t,\mathbf i}
=
\min(\mathcal T_{t,\mathbf i},C_{t,\mathbf i}),
\quad
\delta_{t,\mathbf i}
=
\mathbf 1_{\{\mathcal T_{t,\mathbf i}\leqslant C_{t,\mathbf i}\}},
$
as defined in Section~\ref{sec2}. Figure~\ref{fig:sim_domain} illustrates the spatio-temporal nature of our data.

\subsubsection{Parameter selection and criteria.}
Two complementary criteria are used at distinct stages. Block cross-validation is first used to select the bandwidth \(h\) of the NW estimator or the number of neighbors \(k\) of the \(k\)-NN estimator \citep{hart1994automated,ouyang2006cross}. Once this parameter has been fixed, estimator accuracy is evaluated using the ISE and then compared through its Monte Carlo mean, the MISE. Because a random partition may split neighboring and strongly correlated observations between the training and validation sets, it may lead to an overly optimistic assessment. We therefore use cross-validation based on contiguous blocks, which is appropriate for dependent and spatially autocorrelated data \cite{burman1994cross,schratz2019hyperparameter,rabinowicz2022cross}. Each of the two spatial axes and the temporal axis is divided into two contiguous intervals of equal length, with points lying exactly at the center assigned to the lower interval. The Cartesian product of these partitions thus defines eight spatio-temporal blocks, which are successively left out for validation. The smoothing parameter is selected using a least-squares
		cross-validation criterion, in the spirit of the criteria proposed by
		\cite{patil1993least} and \cite{wang2012hazard} for censored hazard
		function estimation. We adapt it here to the spatio-temporal setting by
		replacing leave-one-out validation with leave-one-block-out validation
		based on contiguous blocks,
	computed up to a model-specific horizon \(\xi\) (a bound beyond
	which the denominator \(1-H\) becomes too unstable). The
	observations are divided into \(B=8\) contiguous blocks
	\(\mathcal I_1,\ldots,\mathcal I_B\). For each block \(b\), the estimator
	\(\widehat\lambda_{\vartheta}^{(-b)}\), with \(\vartheta=h\) for NW or
	\(\vartheta=k\) for \(k\)-NN, is constructed from all observations not
	belonging to \(\mathcal I_b\) and then evaluated on the held-out block. Let
	\(\mathcal G=\{x_1,\ldots,x_G\}\) be a regular grid of \(G=120\) points
	over the domain \([a,\xi]\). For \(\vartheta=h\) or \(k\), we consider
	the criterion
\[
	\mathrm{CV}(\vartheta)
	=
	\frac{1}{B}\sum_{b=1}^{B}
	\left[
	\frac{1}{G}\sum_{g=1}^{G}
	\left\{
	\widehat\lambda_{\vartheta}^{(-b)}(x_g)
	\right\}^{2}
	-
	\frac{2}{|\mathcal I_b|}
	\sum_{j\in\mathcal I_b}
	\frac{\delta_j\,\mathbf 1_{\{a\le X_j\le\xi\}}}{1-\widehat H^{(-b)}(X_j)}\,
	\widehat\lambda_{\vartheta}^{(-b)}(X_j)
	\right].
\]
Here, \(X_j\) and \(\delta_j\) denote respectively the observed
	duration and the uncensoring indicator associated with the \(j\)-th element
	of block \(\mathcal I_b\), corresponding to the variables
	\(X_{t,\mathbf i}\) and \(\delta_{t,\mathbf i}\) defined in
	Section~\ref{sec2}. The index \(j\) is used only to simplify notation,
	each block containing different pairs \((t,\mathbf i)\). The weight
	\(1/\{1-\widehat H^{(-b)}(X_j)\}\) is an inverse censoring-probability
	weight: \(\widehat H^{(-b)}\) is the estimator of \(H\) fitted to the same
	training data as \(\widehat\lambda_{\vartheta}^{(-b)}\) (excluding block
	\(\mathcal I_b\)), and this weighting corrects the bias that censoring
	would otherwise introduce into the cross term. The selected smoothing parameter is defined by
	$
		\widehat{\vartheta}
		=
		\underset{\vartheta\in\Theta}{\arg\min}\,
		\mathrm{CV}(\vartheta),
	$
	where \(\Theta=\mathcal H\) for the Nadaraya--Watson estimator
		and \(\Theta=\mathcal K_n\) for the \(k\)-nearest neighbors estimator.
	The candidate grids are given by
	$
	\mathcal H
	=
	\left\{
	\xi c:\;
	c\in
	\{0.035,\ 0.055,\ 0.080,\ 0.115,\ 0.165,\ 0.230\}
	\right\},
	$
	for the Nadaraya--Watson estimator, and by
	$
	\mathcal K_n
	=
	\left\{
	\max\!\left(
	7,\left\lfloor nc+\frac12\right\rfloor
	\right):
	\;
	c\in
	\{0.035,\ 0.060,\ 0.095,\ 0.145,\ 0.215,\ 0.300\}
	\right\},
	$
	for the \(k\)-nearest neighbors estimator, where
		\(\lfloor x\rfloor\) denotes the integer part of \(x\).
The domain \([a,\xi]\subset\mathcal Z\) is chosen so as to
	exclude the extreme part of the support where the denominator \(1-H\)
	becomes too small and hazard function estimation becomes unstable. We use
	\([0.10,0.80]\) for \(M_1\) and \([0.05,1.20]\) for \(M_2\).
Once \(h\) or \(k\) has been selected, estimator accuracy is
	evaluated on the same grid \(\mathcal G=\{x_1,\ldots,x_G\}\). For
	replication \(r\), we consider the discrete approximation of the normalized
	integrated squared error
$
	\mathrm{ISE}_r
	=
	\dfrac{1}{G}
	\displaystyle\sum_{g=1}^{G}
	\left\{
	\widehat\lambda_r(x_g)-\lambda(x_g)
	\right\}^{2}.
$
Finally, the overall performance of each method is summarized by
	the \textbf{MISE} (\emph{Mean Integrated Squared Error}), defined as the
	empirical mean of the \(\mathrm{ISE}_r\) over the \(250\) independent
	replications:
$
	\mathrm{MISE}
	=
	\dfrac{1}{250}
	\displaystyle\sum_{r=1}^{250}
	\mathrm{ISE}_r.
$
\subsubsection{\texorpdfstring{Main comparison of NW and \(k\)-NN.}{Main comparison of NW and k-NN.}}

\begin{table}[H]
\centering
\scriptsize
\setlength{\tabcolsep}{3.9pt}
\renewcommand{\arraystretch}{0.69}
\caption{MISE of the proposed estimators under the increasing-domain design. Parentheses report the standard deviation of the ISE over 250 replications.}
\label{tab:main_revised}
\begin{tabular}{cccrrr}
\toprule
Model & Cens. & Dependence & \(m^3\) & NW & \(k\)-NN\\
\midrule
\multirow{12}{*}{\(M_1\)}
&20\%&weak&\(5^3\)&0.1539 (0.1059)&0.1010 (0.0936)\\
&&&\(8^3\)&0.0351 (0.0179)&0.0300 (0.0166)\\
&&&\(10^3\)&0.0398 (0.0129)&0.0246 (0.0100)\\
&&strong&\(5^3\)&0.2135 (0.3863)&0.1868 (0.3365)\\
&&&\(8^3\)&0.0806 (0.0660)&0.0785 (0.0660)\\
&&&\(10^3\)&0.0572 (0.0340)&0.0596 (0.0351)\\
&60\%&weak&\(5^3\)&0.1228 (0.1179)&0.1456 (0.1332)\\
&&&\(8^3\)&0.0785 (0.0431)&0.0726 (0.0404)\\
&&&\(10^3\)&0.0581 (0.0281)&0.0577 (0.0264)\\
&&strong&\(5^3\)&0.3491 (0.6439)&0.2582 (0.3525)\\
&&&\(8^3\)&0.1239 (0.1053)&0.0860 (0.0795)\\
&&&\(10^3\)&0.1156 (0.0622)&0.0567 (0.0411)\\
\midrule
\multirow{12}{*}{\(M_2\)}
&20\%&weak&\(5^3\)&0.0314 (0.0311)&0.0347 (0.0370)\\
&&&\(8^3\)&0.0104 (0.0086)&0.0112 (0.0090)\\
&&&\(10^3\)&0.0098 (0.0049)&0.0174 (0.0067)\\
&&strong&\(5^3\)&0.1279 (0.1665)&0.1584 (0.1988)\\
&&&\(8^3\)&0.0373 (0.0474)&0.0434 (0.0508)\\
&&&\(10^3\)&0.0226 (0.0279)&0.0261 (0.0306)\\
&60\%&weak&\(5^3\)&0.0547 (0.0652)&0.0524 (0.0546)\\
&&&\(8^3\)&0.0140 (0.0117)&0.0135 (0.0110)\\
&&&\(10^3\)&0.0084 (0.0065)&0.0061 (0.0051)\\
&&strong&\(5^3\)&0.1753 (0.2556)&0.1671 (0.2672)\\
&&&\(8^3\)&0.0441 (0.0454)&0.0448 (0.0494)\\
&&&\(10^3\)&0.0277 (0.0271)&0.0282 (0.0289)\\
\bottomrule
\end{tabular}
\end{table}

Table~\ref{tab:main_revised} shows an overall improvement in the accuracy of both estimators as the domain size increases. In most scenarios, the MISE decreases as \(m\) increases from \(5\) to \(8\) and then to \(10\), together with a reduction in ISE variability. Strong dependence generally worsens performance relative to weak dependence, especially for small sample sizes. For example, under \(M_2\) with 20\,\% censoring and \(m=5\), the MISE of NW increases from \(0.0314\) to \(0.1279\). The large standard deviations observed in these configurations also reflect greater variability. The effect of the censoring rate is less regular. Under \(M_1\), 60\,\% censoring generally leads to an increase in MISE, particularly in strong-dependence configurations. Under \(M_2\), however, some configurations yield errors comparable to, or even smaller than, those obtained with 20\,\% censoring. Finally, neither estimator dominates uniformly. The \(k\)-NN estimator often performs better for \(M_1\), particularly under heavy censoring, whereas NW generally yields the smallest MISE values for \(M_2\) with 20\,\% censoring. These results confirm the complementarity between the fixed smoothing of NW and the adaptive nature of \(k\)-NN, depending on the hazard shape and the strength of dependence.
\subsubsection{Benchmark estimators.}
We add two benchmark estimators: a smoothed Nelson--Aalen estimator
	\cite{ramlau1983smoothing} and a parametric Weibull estimator
	\cite{weibull1951statistical}. The \(N=\hat{\mathbf n}T\) pooled observations are reindexed as \((X_\ell,\delta_\ell)\), \(\ell=1,\ldots,N\). If \(t_1<\cdots<t_J\) are the distinct failure times, let \(d_j\) denote the number of events at \(t_j\) and \(R_j\) the number at risk immediately before \(t_j\). The smoothed Nelson--Aalen hazard estimator is then
$
\widehat\lambda_{NA,h}(x)
=\dfrac{1}{h}\displaystyle\sum_{j=1}^{J}
K\!\left(\frac{x-t_j}{h}\right)\frac{d_j}{R_j}.
$ The smoothing parameter \(h\) is selected by cross-validation according to the procedure described above.
For the parametric benchmark, the duration \(\mathcal T\) is assumed to follow a Weibull distribution with shape parameter \(a>0\) and scale parameter \(b>0\).
The parameters are estimated from the censored observations using the \texttt{survreg} function from the \texttt{survival} package in \textsf{R}.
The parametric hazard estimator is then
$
	\widehat\lambda_W(x)
	=\dfrac{\widehat a}{\widehat b}
	\left(\dfrac{x}{\widehat b}\right)^{\widehat a-1},
	\quad x>0.
$
Table~\ref{tab:benchmarks} presents the comparison for \(m=8\),
corresponding to the intermediate sample size.

\begin{table}[H]
\centering
\scriptsize
\setlength{\tabcolsep}{2.8pt}
\renewcommand{\arraystretch}{0.94}
\caption{MISE of the proposed methods and benchmark estimators for \(m=8\). Parentheses report the standard deviation of the ISE over 250 replications.}
\label{tab:benchmarks}
\begin{tabular}{ccccccc}
\toprule
Model&Cens.&Dep.&NW&\(k\)-NN&Nelson--Aalen&Weibull\\
\midrule
\(M_1\)&20\%&weak&0.0351 (0.0179)&0.0300 (0.0166)&0.0349 (0.0177)&0.0041 (0.0053)\\
&20\%&strong&0.0806 (0.0660)&0.0785 (0.0660)&0.0800 (0.0649)&0.0376 (0.0550)\\
&60\%&weak&0.0785 (0.0431)&0.0726 (0.0404)&0.0777 (0.0421)&0.0103 (0.0149)\\
&60\%&strong&0.1239 (0.1053)&0.0860 (0.0795)&0.1217 (0.1012)&0.0426 (0.0570)\\
\(M_2\)&20\%&weak&0.0104 (0.0086)&0.0112 (0.0090)&0.0105 (0.0086)&0.0516 (0.0200)\\
&20\%&strong&0.0373 (0.0474)&0.0434 (0.0508)&0.0377 (0.0481)&0.0631 (0.0428)\\
&60\%&weak&0.0140 (0.0117)&0.0135 (0.0110)&0.0142 (0.0125)&0.0200 (0.0152)\\
&60\%&strong&0.0441 (0.0454)&0.0448 (0.0494)&0.0449 (0.0469)&0.0490 (0.0479)\\
\bottomrule
\end{tabular}
\end{table}

Table~\ref{tab:benchmarks} shows that no method dominates uniformly. For \(M_1\), the Weibull model systematically achieves the smallest MISE values, from \(0.0041\) to \(0.0426\), reflecting the good fit of its parametric form to the hazard under consideration. Among the nonparametric methods, \(k\)-NN generally performs best, particularly under 60\,\% censoring and strong dependence, with a MISE of \(0.0860\), compared with \(0.1239\) for NW and \(0.1217\) for Nelson--Aalen. For \(M_2\), the situation is reversed: the nonparametric methods generally outperform Weibull. Under 20\,\% censoring and weak dependence, NW achieves a MISE of \(0.0104\), compared with \(0.0516\) for Weibull. NW and Nelson--Aalen also perform very similarly across all configurations, whereas \(k\)-NN becomes slightly better under 60\,\% censoring and weak dependence. Strong dependence generally increases the MISE and ISE variability for all methods. Overall, these results show that Weibull performs very well when its parametric form is appropriate, as for \(M_1\), but deteriorates when this assumption is unsuitable, as for \(M_2\). The nonparametric estimators therefore appear more flexible, with relatively stable performance across the different scenarios.
\subsubsection{Robustness and violations of assumptions.}
We perform 250 replications for \(m=8\) and 40\,\% censoring in eight scenarios covering, in particular, nonseparability, nonstationarity, anisotropy, irregular hazards, regional heterogeneity, and informative censoring.
These three covariance structures are studied under \(M_2\). The separable structure is
	\[
	\operatorname{Cov}\{Z(s,t),Z(s',t')\}
	=
	\exp\{-0.70\|s-s'\|\}\exp\{-0.70|t-t'|\}.
	\]
The nonseparable structure is
\[
	\operatorname{Cov}\{Z(s,t),Z(s',t')\}
	=
	\exp\left\{
	-\left\|
	\left(
		\frac{s_1-s_1'}{2.2},
		\frac{s_2-s_2'}{2.2},
		\frac{t-t'}{2.0}
	\right)
	\right\|
	\right\}.
\]
The anisotropic and nonstationary structure is
\[
	\varphi(s,t)=
	\left(
		\frac{s_1}{1.3}+0.035\,s_1^{2},
		\ \frac{s_2}{3.0},
		\ \frac{t}{2.0}+0.025\,t^{2}
	\right).
\]

\[
	\operatorname{Cov}\{Z(s,t),Z(s',t')\}
	=
	\exp\bigl\{-\|\varphi(s,t)-\varphi(s',t')\|\bigr\}.
\]
The peaked and discontinuous hazard scenario is
\[
	\lambda_{\mathrm{pic}}(x)=
	\begin{cases}
		0.4, & \text{if } x<0.7,\\
		2.2, & \text{if } 0.7\leqslant x\leqslant 1.1,\\
		0.6, & \text{if } x>1.1.
	\end{cases}
\]
The heterogeneous-hazard scenario is
\[
	\lambda_{\mathrm{het}}(x)=
	\frac{
		\displaystyle
		\sum_{\theta\in\{0.60,\,1.40\}}
		0.5\,\exp(-\theta x^{0.7})
		\times 0.7\,\theta\,x^{-0.3}
	}{
		\displaystyle
		\sum_{\theta\in\{0.60,\,1.40\}}
		0.5\,\exp(-\theta x^{0.7})
	}.
\]
Two additional models complete the hazard shapes under study: \(M_3\), corresponding to an increasing Gompertz-type hazard, and \(M_4\), exhibiting a bathtub shape,
\[
	\lambda_3(x)=0.35\exp(0.70x),
	\quad
	\lambda_4(x)=0.30+1.20\exp(-3x)+0.18x.
\]
\begin{table}[H]
\centering
\scriptsize
\setlength{\tabcolsep}{3pt}
\renewcommand{\arraystretch}{0.95}
\caption{MISE in the robustness study (\(m=8\), 40\,\% censoring, 250 replications). Parentheses report the standard deviation of the ISE.}
\label{tab:robustness}
\begin{tabular}{lrrrr}
\toprule
Scenario&NW&\(k\)-NN&Nelson--Aalen&Weibull\\
\midrule
Stationary separable&0.040 (0.044)&0.043 (0.048)&0.041 (0.047)&0.052 (0.044)\\
Nonseparable covariance&0.147 (0.260)&0.156 (0.274)&0.150 (0.268)&0.158 (0.296)\\
Anisotropy and nonstationarity&0.062 (0.074)&0.066 (0.078)&0.064 (0.078)&0.061 (0.055)\\
Peaked and discontinuous hazard&0.151 (0.128)&0.174 (0.153)&0.148 (0.123)&0.573 (0.062)\\
Regionally heterogeneous hazards&0.110 (0.062)&0.084 (0.056)&0.110 (0.061)&0.030 (0.035)\\
Informative censoring&0.069 (0.056)&0.068 (0.055)&0.067 (0.055)&0.136 (0.065)\\
Increasing Gompertz hazard&0.035 (0.044)&0.032 (0.034)&0.041 (0.058)&0.036 (0.030)\\
Bathtub-shaped hazard&0.109 (0.131)&0.034 (0.033)&0.075 (0.070)&0.036 (0.025)\\
\bottomrule
\end{tabular}
\end{table}

\noindent
Table~\ref{tab:robustness} shows that no method dominates uniformly. In the stationary separable scenario, NW has the smallest MISE (\(0.040\)), followed by Nelson--Aalen (\(0.041\)) and \(k\)-NN (\(0.043\)). Under nonseparable covariance, the errors increase sharply for all methods, with a slight advantage for NW (\(0.147\)), while the large standard deviations indicate substantial variability. In the presence of anisotropy and nonstationarity, the four methods remain close, with Weibull achieving the smallest MISE (\(0.061\)). For a peaked and discontinuous hazard, Nelson--Aalen (\(0.148\)) and NW (\(0.151\)) perform best, whereas Weibull deteriorates sharply (\(0.573\)), illustrating the limitations of a parametric model for an irregular hazard function. For regionally heterogeneous hazards, Weibull achieves the smallest MISE (\(0.030\)), while \(k\)-NN (\(0.084\)) outperforms NW and Nelson--Aalen (\(0.110\)). Under informative censoring, the three nonparametric methods remain very close, with MISE values between \(0.067\) and \(0.069\), compared with \(0.136\) for Weibull. Finally, \(k\)-NN performs best for the increasing Gompertz hazard (\(0.032\)) and especially for the bathtub-shaped hazard (\(0.034\)), compared with \(0.109\) for NW in the latter case. These results suggest that its local adaptation can be advantageous for complex hazard shapes. Overall, NW remains stable near the reference setting, \(k\)-NN becomes particularly competitive in some nonlinear or heterogeneous configurations, Nelson--Aalen provides a robust benchmark, whereas Weibull depends more strongly on how well its parametric form matches the true hazard.

\subsubsection{Sensitivity to the choice of kernel.}
We compare the Epanechnikov, biweight, and triweight kernels, all of which satisfy \hyperref[H2]{H2}, in four representative scenarios: a hazard under weak dependence and a hazard under strong dependence (both under model \(M_2\)), a peaked and discontinuous hazard, and a bathtub-shaped hazard (\(M_4\)). For each kernel, \(h\) and \(k\) are selected as above.

\begin{table}[H]
\centering
\scriptsize
\setlength{\tabcolsep}{4.2pt}
\renewcommand{\arraystretch}{0.70}
\caption{Sensitivity of the MISE to the choice of kernel for \(m=8\). Parentheses report the standard deviation of the ISE over 250 replications.}
\label{tab:kernel_sensitivity}
\begin{tabular}{llrrr}
\toprule
Scenario&Estimator&Epanechnikov&Biweight&Triweight\\
\midrule
\multirow{2}{*}{Weak dependence, 40\,\% censoring}&NW&0.0100 (0.0082)&0.0114 (0.0086)&0.0128 (0.0090)\\
&\(k\)-NN&0.0137 (0.0092)&0.0118 (0.0087)&0.0132 (0.0091)\\
\multirow{2}{*}{Strong dependence, 60\,\% censoring}&NW&0.0463 (0.0543)&0.0495 (0.0595)&0.0522 (0.0633)\\
&\(k\)-NN&0.0446 (0.0516)&0.0472 (0.0535)&0.0497 (0.0556)\\
\multirow{2}{*}{Peak and discontinuities}&NW&0.1349 (0.0988)&0.1382 (0.1039)&0.1321 (0.0966)\\
&\(k\)-NN&0.1508 (0.1273)&0.1474 (0.1231)&0.1426 (0.1186)\\
\multirow{2}{*}{Bathtub-shaped hazard}&NW&0.0860 (0.0793)&0.0748 (0.0670)&0.0816 (0.0747)\\
&\(k\)-NN&0.0317 (0.0227)&0.0347 (0.0237)&0.0377 (0.0249)\\
\bottomrule
\end{tabular}
\end{table}

Table~\ref{tab:kernel_sensitivity} shows that the choice of kernel affects performance without changing the general trends. In the weak-dependence scenario with 40\,\% censoring, Epanechnikov performs best for NW (\(0.0100\)), whereas the biweight kernel yields the smallest MISE for \(k\)-NN (\(0.0118\)).

Under strong dependence and 60\,\% censoring, Epanechnikov performs best for both estimators, with MISE values of \(0.0463\) for NW and \(0.0446\) for \(k\)-NN. However, the differences between kernels remain moderate compared with the effects of dependence and censoring.

For the peaked and discontinuous hazard, the triweight kernel yields the smallest MISE values, namely \(0.1321\) for NW and \(0.1426\) for \(k\)-NN. In the bathtub scenario, biweight is preferable for NW (\(0.0748\)), whereas Epanechnikov remains clearly better for \(k\)-NN (\(0.0317\)).

Overall, no kernel dominates in every configuration. Epanechnikov nevertheless appears to be the most consistent choice, particularly for \(k\)-NN, whereas biweight or triweight may be advantageous for some more complex hazard shapes. Sensitivity to the kernel remains smaller than sensitivity to dependence, censoring, or hazard shape.
\subsubsection{Separate spatial and temporal dependence.}
To distinguish the spatial and temporal effects, we consider \(M_2\) with \(m=8\) and 40\,\% censoring, combining weak and strong levels of the parameters \(a\) and \(b\). The parameters \(h\) and \(k\) are selected as above.
\begin{table}[H]
\centering
\scriptsize
\setlength{\tabcolsep}{4pt}
\renewcommand{\arraystretch}{0.7}
\caption{MISE when spatial and temporal dependence vary separately (\(M_2\), \(m=8\), 40\,\% censoring, 250 replications). Parentheses report the Monte Carlo errors of the MISE.}
\label{tab:separate_dependence}
\begin{tabular}{ccrrr}
\toprule
Spatial dep.&Temporal dep.&NW&\(k\)-NN&Nelson--Aalen\\
\midrule
weak&weak&0.0169 (0.0007)&0.0188 (0.0007)&0.0172 (0.0007)\\
strong&weak&0.0261 (0.0024)&0.0273 (0.0026)&0.0265 (0.0025)\\
weak&strong&0.0122 (0.0007)&0.0123 (0.0007)&0.0121 (0.0007)\\
strong&strong&0.0421 (0.0041)&0.0492 (0.0047)&0.0471 (0.0044)\\
\bottomrule
\end{tabular}
\end{table}

Table~\ref{tab:separate_dependence} shows that the effects of spatial and temporal dependence are not symmetric. When temporal dependence remains weak, moving from weak to strong spatial dependence substantially increases the MISE: for NW, it rises from \(0.0169\) to \(0.0261\), and for \(k\)-NN from \(0.0188\) to \(0.0273\). Strong spatial dependence therefore appears to be an important source of loss of accuracy. Conversely, when spatial dependence is weak, strong temporal dependence does not degrade performance in this experiment and even leads to smaller MISE values, around \(0.012\) for all three methods. However, this effect should be interpreted as an empirical result specific to this configuration and not as a general improvement induced by temporal dependence. When both forms of dependence are strong, the errors increase sharply for all methods, with MISE values of \(0.0421\) for NW, \(0.0492\) for \(k\)-NN, and \(0.0471\) for Nelson--Aalen. These results therefore suggest that the deterioration is mainly associated with strong spatial dependence and becomes more pronounced when it is combined with strong temporal dependence. NW then has the smallest MISE among the three estimators.

\subsection{\texorpdfstring{Applications to Real Data}{Applications to Real Data}}

We consider two complementary applications: the more controlled and homogeneous CGD trial and the richer but more heterogeneous spatio-temporal NYC~311 data. They illustrate two distinct contexts in which the proposed estimators can be applied.

\subsubsection{\texorpdfstring{First application: multicenter clinical trial with observed right censoring}{First application: multicenter clinical trial with observed right censoring}}

\begin{figure}[H]
		\centering
		\includegraphics[width=0.8\textwidth]{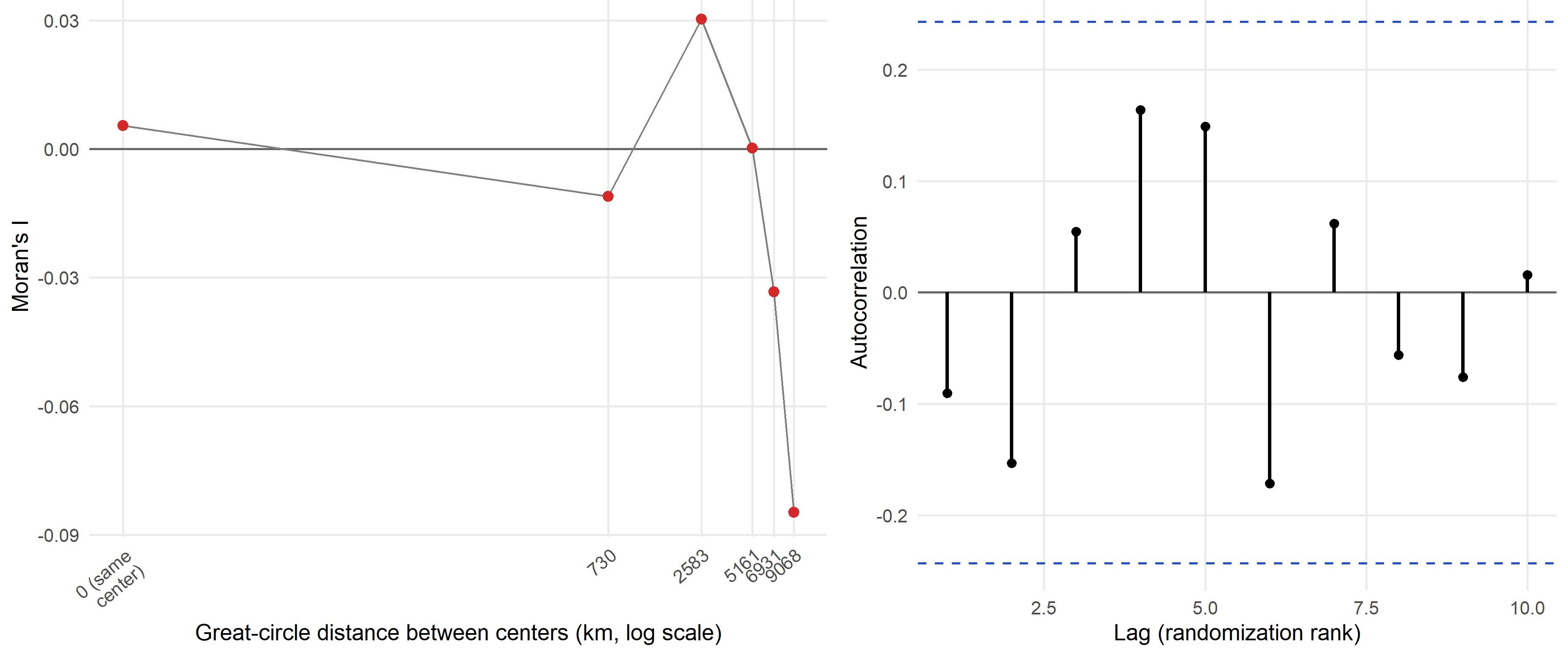}
		\caption{Dependence structure of the placebo cohort: correlation between centers (left) and temporal autocorrelation (right).}
		\label{fig:cgd_dependance}
	\end{figure}

We use data from the randomized, double-blind clinical trial on chronic granulomatous disease (CGD), conducted among 128 patients at 13 centers \cite{cgd1991}. The outcome of interest is the time to the first serious infection; in the absence of an event, the observation is right-censored at the end of follow-up. The data are available in the \texttt{cgd0} data set of the \texttt{survival} package in \textsf{R} \cite{therneau2026survival}. The analysis is restricted to the placebo group in order to remove the treatment effect, whose efficacy was established in the original trial \cite{cgd1991}. Among the 65 placebo patients, 30 experienced at least one serious infection and 35 were censored, corresponding to a censoring rate of 53.8\,\%. Each patient \(j=1,\ldots,65\) is characterized by the enrollment center \(i_j\) and randomization date \(t_j\), between June 7 and December 29, 1988. The resulting design is irregular, as several patients may belong to the same center or even share the same randomization date. Letting \(\mathcal T_j\) denote the time to the first serious infection and \(C_j\) the censoring time, we observe
	\[
	X_j=\min(\mathcal T_j,C_j),
	\qquad
	\delta_j=\mathbf 1_{\{\mathcal T_j\leqslant  C_j\}}.
	\]
	Thus, \(\delta_j=1\) corresponds to an observed infection and \(\delta_j=0\) to censoring.
The smoothing parameters are selected by five-block cross-validation, with the blocks constructed from center--month groups and then allocated so as to balance their sizes.
The proposed estimators are compared with two standard benchmarks, one nonparametric and the other parametric. Figure~\ref{fig:cgd_dependance} highlights the dependence structures associated with center and randomization date.
Moran's indices remain close to zero for all distance classes, while the temporal autocorrelations remain within the confidence bands. These results suggest overall weak spatio-temporal dependence in the observed durations.

	\begin{figure}[H]
		\centering
		\includegraphics[width=0.5\textwidth]{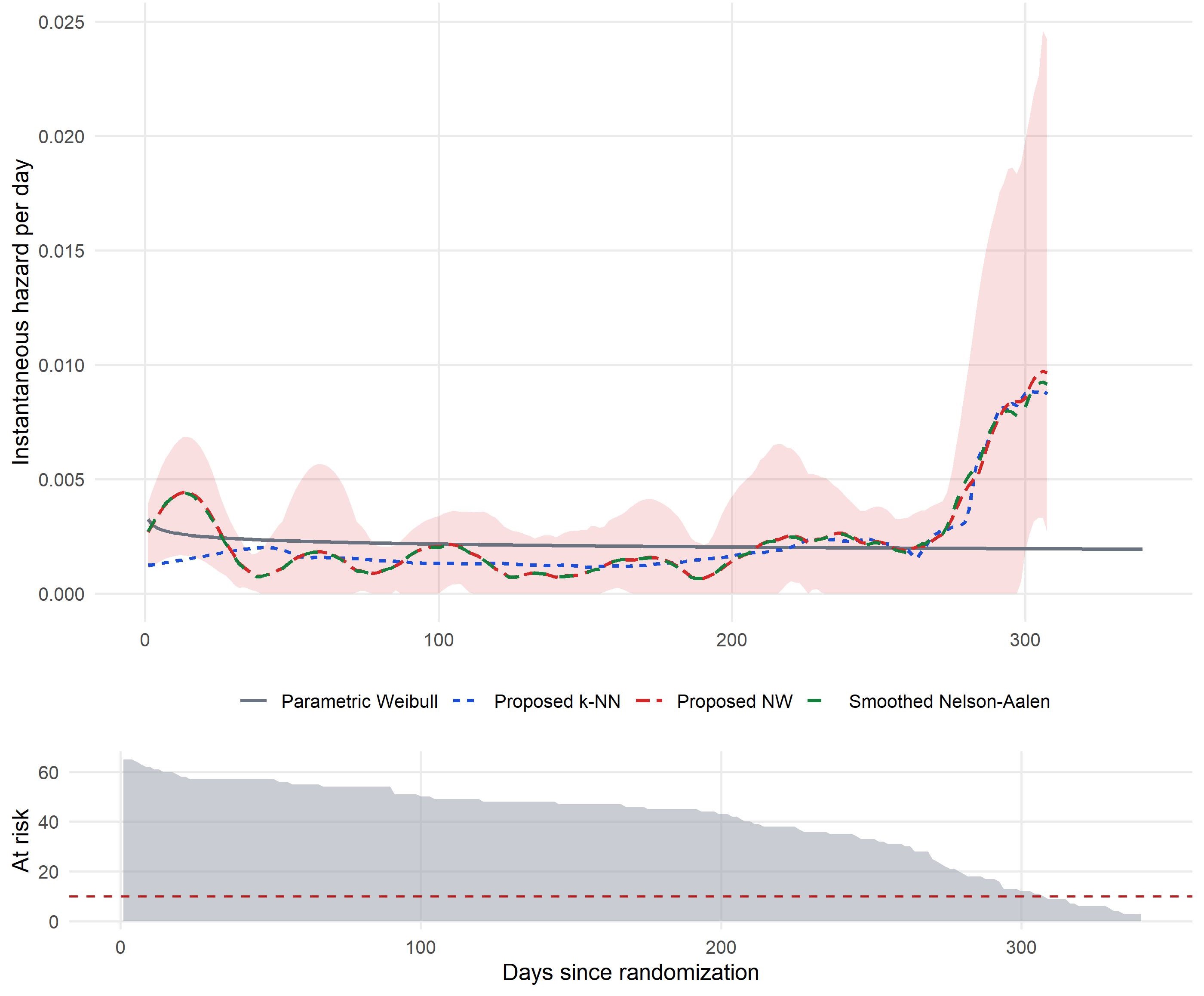}
		\caption{Estimates of the marginal hazard of a first serious infection using the four methods. The lower panel shows the number of patients still at risk.}
		\label{fig:app_hazard}
	\end{figure}

The NW, \(k\)-NN, and smoothed Nelson--Aalen curves are truncated after 307 days, when the number at risk becomes insufficient. Over this interval, NW and Nelson--Aalen are almost indistinguishable and reveal an early peak, while \(k\)-NN follows a similar pattern. The more rigid Weibull model reproduces these variations less accurately.
\subsubsection{\texorpdfstring{Second application: geolocated and censored NYC~311 requests}{Second application: geolocated and censored NYC 311 requests}}

\begin{figure}[H]
		\centering
		\includegraphics[width=0.8\textwidth]{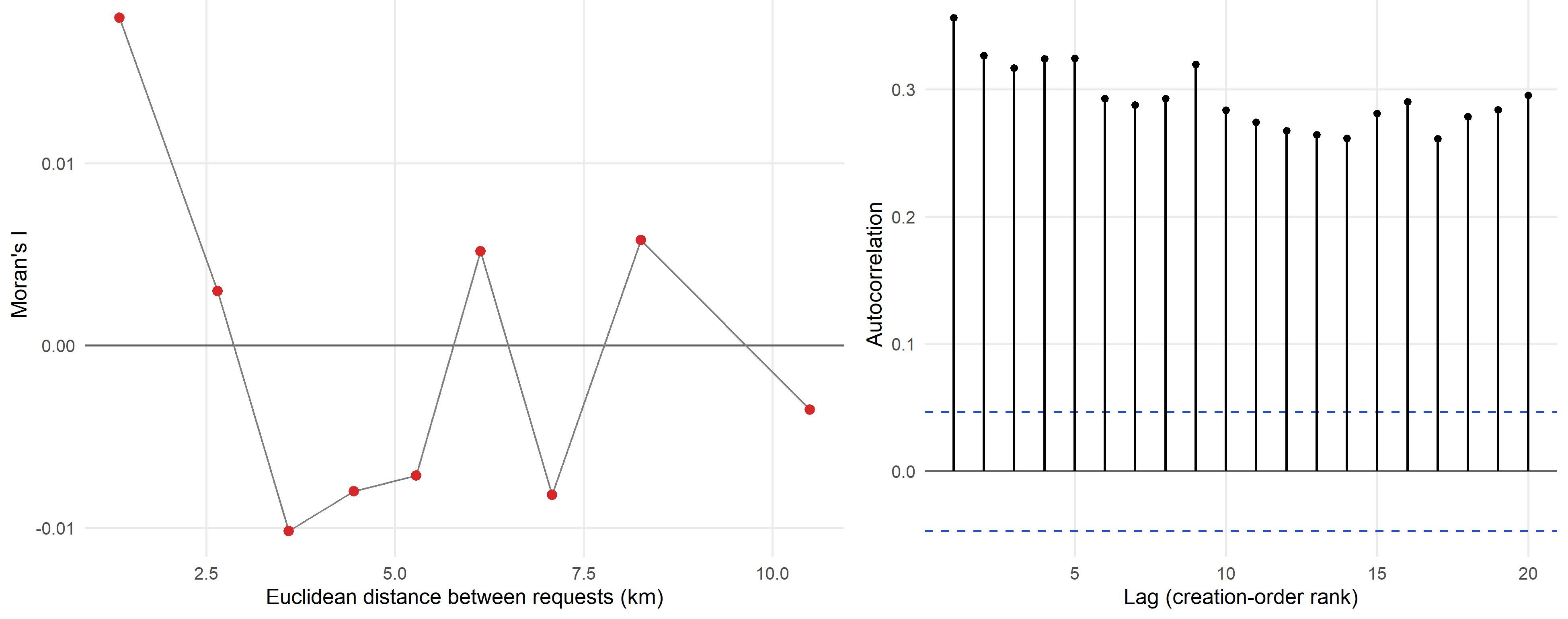}
		\caption{Spatio-temporal dependence structure of the NYC~311 data: spatial correlogram (left) and temporal autocorrelation function (right).}
		\label{fig:nyc_dependance}
	\end{figure}

NYC~311 is New York City's non-emergency municipal service. We use the official \emph{311 Service Requests from 2020 to Present} data set \cite{nyc3112026}, frozen on August 14, 2026 together with its metadata and SHA-256 checksum. The data and their dictionary are freely available on the \href{https://data.cityofnewyork.us/d/erm2-nwe9}{NYC Open Data platform}.
The study concerns requests created in Brooklyn in July 2026, classified as \emph{Damaged Tree} with the detail \emph{Branch or Limb Has Fallen Down}, having complete coordinates and a status of \emph{Closed} or \emph{In Progress}. At the cutoff date of August 13, 2026, 1,127 of the 1,746 requests are closed and 619 are right-censored, corresponding to a censoring rate of 35.5\,\%. Each request \(j=1,\ldots,1746\) is associated with its creation date \(t_j\) and geographic position \(\mathbf s_j=(\mathrm{longitude}_j,\mathrm{latitude}_j)\). Letting \(\mathcal T_j\) denote the time to closure and \(C_j=t^\star-t_j\) the follow-up available up to the cutoff date \(t^\star\), we observe
\[
	X_j=\min(\mathcal T_j,C_j),
	\qquad
	\delta_j=\mathbf 1_{\{\mathcal T_j\leqslant C_j\}}.
\]
The median observed duration is 2.58 days, compared with 29.45 days of follow-up for the censored requests. The smoothing parameters are selected by cross-validation using five spatio-temporal blocks. Figure~\ref{fig:nyc_dependance} illustrates the spatial and temporal dependence structures observed in the data.
The spatial correlogram indicates weak dependence, mainly at short distances, with Moran's indices remaining generally close to zero. By contrast, the temporal autocorrelations are positive, large, and persistent across all lags considered. Temporal dependence therefore appears substantially stronger than spatial dependence.
	\begin{figure}[H]
		\centering
		\includegraphics[width=0.48\textwidth]{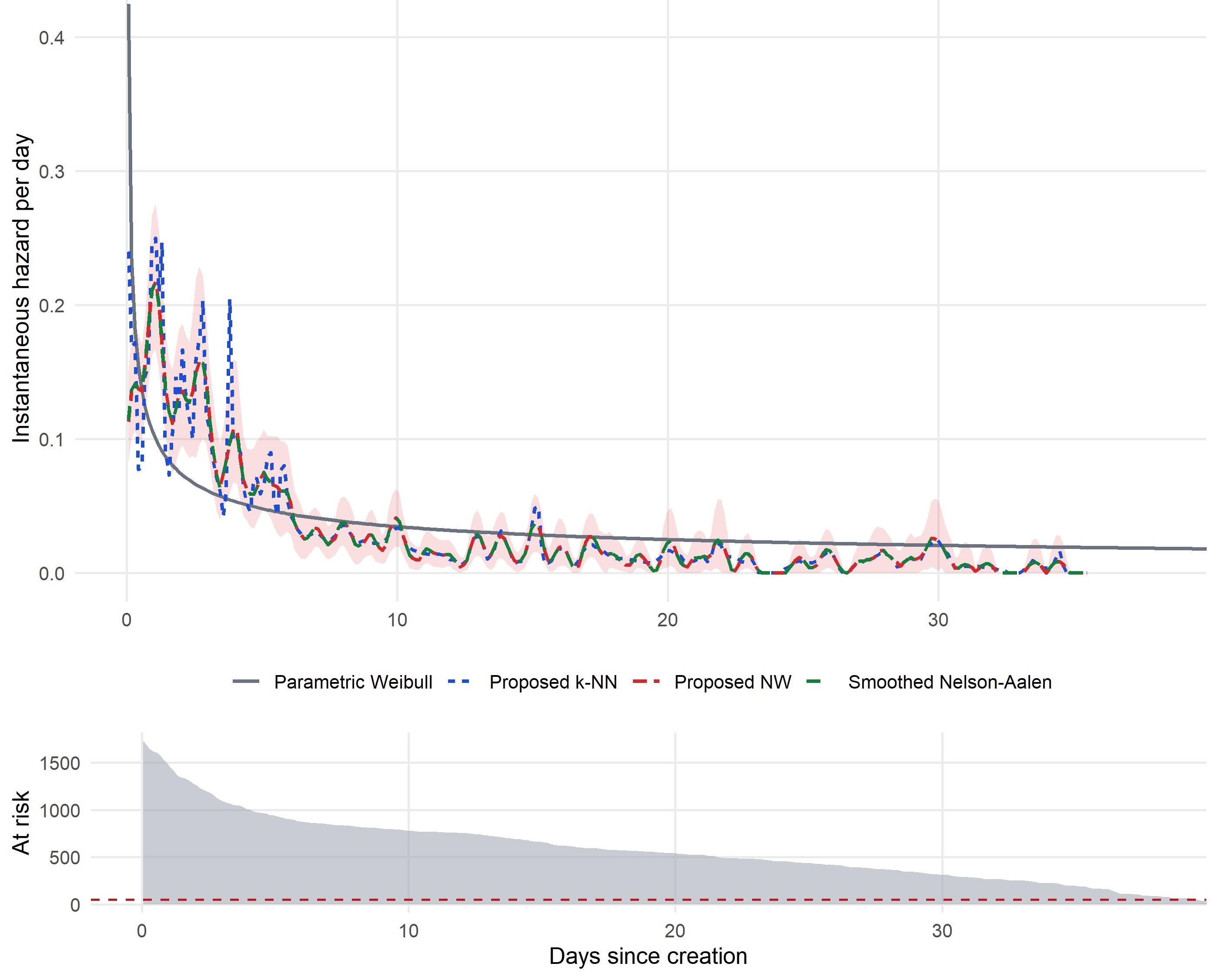}
		\caption{Estimates of the marginal closure hazard for NYC~311 requests using the four methods. The lower panel shows the numbers at risk.}
		\label{fig:nyc_hazard}
	\end{figure}

Figure~\ref{fig:nyc_hazard} shows a high closure hazard at the beginning of follow-up, followed by a sharp decline after six days. NW and smoothed Nelson-Aalen are almost indistinguishable, while \(k\)-NN exhibits more pronounced local variation and Weibull captures the overall trend. Beyond six days, the nonparametric estimators become very close.

 Both figures show 95\,\% pointwise NW bootstrap bands and risk-set panels; uncertainty increases as risk sets shrink. The nonparametric curves in Figure~\ref{fig:nyc_hazard} are truncated when the observed survival falls below \(0.10\). In summary, the two applications confirm the relevance of the proposed estimators. The \(k\)-NN estimator adapts better to pronounced local variation, as in NYC~311, whereas NW provides a smoother and more stable estimate, particularly in the CGD application.
\section{Conclusion}
\label{sec5}

\noindent
This article studied the nonparametric estimation of a common marginal hazard function from right-censored data exhibiting spatial and temporal dependence. The target \(\lambda(x)\) depends only on the duration \(x\), while the spatial and temporal indices enter the dependence structure. Two estimators, based respectively on Nadaraya--Watson and the \(k\)-nearest neighbors method, were proposed. Under the stated assumptions, their pointwise and uniform consistency is established, with the same order of asymptotic convergence. The simulation study, conducted under an increasing spatio-temporal domain, confirms the importance of effective sample size and the value of block cross-validation for limiting the effects of dependence. The results provide no evidence of systematic superiority of one method over the others. The Weibull model is naturally favored when correctly specified, whereas NW appears to be the most stable nonparametric estimator overall. The \(k\)-NN estimator stands out more clearly when the hazard function exhibits pronounced local variation. Stronger dependence generally degrades performance, and informative censoring increases bias as the dependence between the event and censoring becomes stronger. The robustness experiments also show that the theoretical results do not automatically extend to nonseparable, nonstationary, or heterogeneous settings. By contrast, the choice of kernel does not change the general conclusions regarding the relative behavior of the estimators. The applications to the CGD and NYC~311 data illustrate these results in two distinct right-censoring settings. The diagnostics indicate weak residual dependence for CGD, whereas the NYC~311 data exhibit short-range spatial dependence and, above all, strong temporal persistence. In the latter case, \(k\)-NN captures local hazard variation more clearly, while NW and smoothed Nelson--Aalen provide smoother estimates and remain particularly close in the CGD application. Finally, the assumption of a common marginal hazard \(\lambda(x)\) presumes a common marginal distribution despite spatio-temporal dependence. In the presence of spatial, temporal, or individual heterogeneity, a natural extension would be to consider a conditional hazard
	$
	\lambda(x\mid \mathbf{s},t,\mathbf{z}),
	$
	explicitly incorporating location, time, and possible covariates.

		\section{Proofs}\label{sec6}
	In this section, we present the proofs of the main results. The technical lemmas supporting these proofs are stated and proved in the subsection below. The proofs of Theorems~\ref{th3} and~\ref{th4}, which are analogous to those of Theorems~\ref{th1} and~\ref{th2}, are omitted.
	 Since weak $\tau$-dependence implies weak $\eta$-dependence (see \cite{dedecker04}), we adopt the latter framework to derive variance and covariance bounds, following the methodology of \cite{amiri18}. Consequently, there exist functions $\eta_1(\cdot)$ and $\eta_2(\cdot)$ satisfying $\eta_1(0) = \eta_2(0) = 1$, such that for all $z, w \in \mathbb{R}$,
	\[
	\eta(z, w) = \eta_1(z)\,\eta_2(w),
	\]
	where
	\[
	\eta_1(z) = O(z^{-\theta}) \quad \text{and} \quad
	\eta_2(t_{\mathbf{j}}) = O(t_{\mathbf{j}}^{-\theta}) \quad \text{for some} \quad \theta > 2,
	\]
	with $t_{\mathbf{j}} = \|\mathbf{i} - \mathbf{j}\|$, for $\mathbf{i}, \mathbf{j} \in \Lambda_{\mathbf{n}}$, and $z = |t - s|$ denoting the temporal lag between two observations.

\subsection{Technical Lemmas}\label{ssec:technical_lemmas}

\begin{lemma}[]\label{lm1}
		Let \(A_{t,\mathbf{i}}\) be an \(\eta\)-weakly dependent spatio-temporal process with bounded density \(l\).
		Under assumptions \(\hyperref[H1]{H1}\), \(\hyperref[H2]{H2}\), and \(\hyperref[H6]{H6}\), let \(p\) be the function defined by
		$
		p(t) = K\left(\dfrac{x - t}{\widetilde{h}_{T,\mathbf{n}}}\right).
		$
		Then,
		\[
		\mathrm{Var}\left(\sum_{\mathbf{i} \in \Lambda_{\mathbf{n}}}\sum_{t=1}^{T} p(A_{t,\mathbf{i}})\right)
		= O\!\left(\widehat{\mathbf{n}}\, T\, \widetilde{h}_{T,\mathbf{n}}\right).
		\]
	\end{lemma}
	
	\begin{proof}[Proof of Lemma~\ref{lm1}]
	For notational simplicity, set
	$
	K_{t,\mathbf{i}} := p(A_{t,\mathbf{i}})
	= K\left(\frac{x-A_{t,\mathbf{i}}}{\widetilde{h}_{T,\mathbf{n}}}\right).
	$
	\begin{eqnarray*}
		Var\left(\displaystyle\sum_{\textbf{i} \in \Lambda_{\mathbf{n}}}^{}\sum_{t=1}^{T}p(A_{t,\textbf{i}})\right) &=& \displaystyle\sum_{\textbf{i} \in \Lambda_{\mathbf{n}}}^{}\sum_{t=1}^{T}Var\left(p(A_{t,\textbf{i}})\right) + \displaystyle\sum_{\textbf{i} \in \Lambda_{\mathbf{n}}}^{}\sum_{s \ne t}^{T}\text{Cov}(p(A_{t,\textbf{i}}),p(A_{s,\textbf{i}}))\\
		&+& \displaystyle\sum_{\substack{\textbf{i},\textbf{j} \in \Lambda_{\mathbf{n}} \\ \textbf{i} \ne \textbf{j}}}
		^{}\sum_{t=1}^{T}\text{Cov}(p(A_{t,\textbf{i}}),p(A_{t,\textbf{j}}))
		+ \displaystyle\sum_{\substack{\textbf{i},\textbf{j} \in \Lambda_{\mathbf{n}} \\ \textbf{i} \ne \textbf{j}}}
		^{}\sum_{s \ne t}^{T}\text{Cov}(p(A_{t,\textbf{i}}),p(A_{s,\textbf{j}}))\\
		&=& \mathcal{I}_1 + \mathcal{I}_2 + \mathcal{I}_3 + \mathcal{I}_4.
	\end{eqnarray*}
	
	For \(\mathcal{I}_1\), we have:
	$
	Var\left(K_{t,\textbf{i}}\right)
	\leqslant \widetilde{h}_{T,\textbf{n}}\displaystyle\int_{D}K^2\left(u\right) l\left(x - u\widetilde{h}_{T,\textbf{n}}\right) du \hspace{0.2cm}
	= O(\widetilde{h}_{T,\textbf{n}})$ because \( K \) has compact support \( D \). Thus, \(\mathcal{I}_1 =  O\left(\widehat{\textbf{n}}T\widetilde{h}_{T,\textbf{n}}\right)\).
	For \(\mathcal{I}_2\), we have:
	\begin{eqnarray*}
		\mathcal{I}_2 &=& \displaystyle\sum_{\textbf{i} \in \Lambda_{\mathbf{n}}}^{}\sum_{s \ne t}^{T}\text{Cov}(p(A_{t,\textbf{i}}),p(A_{s,\textbf{i}}))
		= 2\displaystyle\sum_{\textbf{i} \in \Lambda_{\mathbf{n}}}^{}\sum_{s< t}^{T}\text{Cov}(p(A_{t,\textbf{i}}),p(A_{s,\textbf{i}}))\\
		&=& 2\displaystyle\sum_{\textbf{i} \in \Lambda_{\mathbf{n}}}^{}\sum_{0 < s-t  \leqslant  a_{t,n}}^{T}\text{Cov}(p(A_{t,\textbf{i}}),p(A_{s,\textbf{i}})) +  2\displaystyle\sum_{\textbf{i} \in \Lambda_{\mathbf{n}}}^{}\sum_{s - t > a_{t,n}}^{T}\text{Cov}(p(A_{t,\textbf{i}}),p(A_{s,\textbf{i}}))\\
		&=& \mathcal{J}_1 + \mathcal{J}_2.
	\end{eqnarray*}
	For \(\mathcal{J}_1\), we have:
\[
\begin{aligned}
	\operatorname{Cov}\!\left(K_{t,\mathbf{i}},K_{s,\mathbf{i}}\right)
	&= \mathbb{E}\!\left[K_{t,\mathbf{i}}K_{s,\mathbf{i}}\right]
	- \mathbb{E}\!\left[K_{t,\mathbf{i}}\right]
	\mathbb{E}\!\left[K_{s,\mathbf{i}}\right] \\[0.5em]
	&= \int_D \int_D
	K\!\left(\frac{x-t_1}{\widetilde{h}_{T,\mathbf{n}}}\right)
	K\!\left(\frac{x-t_2}{\widetilde{h}_{T,\mathbf{n}}}\right)
	l_{(X_{t,\mathbf{i}},X_{s,\mathbf{i}})}(t_1,t_2)\,dt_1dt_2 \\[0.5em]
	&\quad
	- \int_D \int_D
	K\!\left(\frac{x-t_1}{\widetilde{h}_{T,\mathbf{n}}}\right)
	K\!\left(\frac{x-t_2}{\widetilde{h}_{T,\mathbf{n}}}\right)
	l(t_1)l(t_2)\,dt_1dt_2 .
\end{aligned}
\]
\[
\begin{aligned}
	\left|
	\operatorname{Cov}\!\left(K_{t,\mathbf{i}},K_{s,\mathbf{i}}\right)
	\right|
	&\leqslant \widetilde{h}_{T,\mathbf{n}}^{2}
	\int_D \int_D
	\left|K(u_1)K(u_2)\right| \\[0.3em]
	&\qquad \times
	\Bigl|
	l_{(X_{t,\mathbf{i}},X_{s,\mathbf{i}})}
	\bigl(x-u_1\widetilde{h}_{T,\mathbf{n}},
	x-u_2\widetilde{h}_{T,\mathbf{n}}\bigr) \\[0.3em]
	&\qquad\qquad
	- l\bigl(x-u_1\widetilde{h}_{T,\mathbf{n}}\bigr)
	l\bigl(x-u_2\widetilde{h}_{T,\mathbf{n}}\bigr)
	\Bigr|\,du_1du_2 \\[0.5em]
	&\leqslant \widetilde{h}_{T,\mathbf{n}}^{2}
	\sup_{(t,\mathbf{i})\neq (s,\mathbf{i})}
	\left\|
	l_{(X_{t,\mathbf{i}},X_{s,\mathbf{i}})} - l\otimes l
	\right\|_{\infty}
	\int_D \int_D
	\left|K(u_1)K(u_2)\right|\,du_1du_2 \\[0.5em]
	&= O\!\left(\widetilde{h}_{T,\mathbf{n}}^{2}\right).
\end{aligned}
\]
	Setting \( k= s- t \), we obtain:
	\[
	\begin{aligned}
		\mathcal{J}_1
		&= 2\sum_{\mathbf{i}\in\Lambda_{\mathbf{n}}}\sum_{s=1}^{T}\sum_{k=1}^{a_{t,n}}
		\operatorname{Cov}\!\bigl(p(A_{t,\mathbf{i}}),p(A_{s,\mathbf{i}})\bigr) \\[0.4em]
		&\leqslant 2\sum_{\mathbf{i}\in\Lambda_{\mathbf{n}}}\sum_{s=1}^{T}\sum_{k=1}^{a_{t,n}}
		\left|\operatorname{Cov}\!\bigl(p(A_{t,\mathbf{i}}),p(A_{s,\mathbf{i}})\bigr)\right| 
		= O\!\left(\widehat{\mathbf{n}}\,T\,a_{t,n}\,\widetilde{h}_{T,\mathbf{n}}^{2}\right).
	\end{aligned}
	\]
	For \( \mathcal{J}_2 \), under hypothesis \hyperref[H2]{H2} and by \cite[Corollary 3.3]{hwang13}, we have:
	\begin{eqnarray*}
		\mathcal{J}_2 &=& 2\displaystyle\sum_{\textbf{i} \in \Lambda_{\mathbf{n}}}^{}\sum_{s = 1}^{T}\sum_{k > a_{t,n}}^{}\text{Cov}(p(A_{t,\textbf{i}}),p(A_{s,\textbf{i}}))
		\leqslant  2\displaystyle\sum_{\textbf{i} \in \Lambda_{\mathbf{n}}}^{}\sum_{s = 1}^{T}\sum_{k > a_{t,n}}^{}\left|\text{Cov}(p(A_{t,\textbf{i}}),p(A_{s,\textbf{i}}))\right|\\
		&\leqslant & 2C\widehat{\textbf{n}}T\left\Vert p(A_{t,\textbf{i}})\right\Vert_4\left\Vert f(A_{s,\textbf{i}})\right\Vert_4\sum_{k > a_{t,n}}^{}(\eta_1(k))^{1/2}\\
		&\leqslant & 2C\widehat{\textbf{n}}T\left\Vert p(A_{t,\textbf{i}})\right\Vert_4\left\Vert p(A_{\tau,\textbf{i}})\right\Vert_4\int_{ a_{t,n}}^{\infty}k^{-\frac{\theta}{2}}dk
		= O\left(\widehat{\textbf{n}} Ta^{1-\frac{\theta}{2}}_{t,n}\widetilde{h}^{1/2}_{T,\textbf{n}}\right).
	\end{eqnarray*}
	Setting \( a_{t,n} = \left\lfloor(\widetilde{h}_{T,\textbf{n}})^{-w} \right\rfloor\) with \( w > 1/(\theta -2) \), we obtain: $\mathcal{I}_2 = \mathcal{J}_1 + \mathcal{J}_2 = O\left(\widehat{\textbf{n}}T\widetilde{h}_{T,\textbf{n}}\right)$. Let us now determine $\mathcal{I}_3$. We have:
	\begin{eqnarray*}
		\mathcal{I}_3 &=& \sum_{\substack{\textbf{i},\textbf{j} \in \Lambda_{\mathbf{n}} \\ \textbf{i} \ne \textbf{j}}}^{}\sum_{t=1}^{T}\text{Cov}(p(A_{t,\textbf{i}}),p(A_{t,\textbf{j}}))\\
		&=&\sum_{\substack{\textbf{i},\textbf{j} \in \Lambda_{\mathbf{n}} \\ \Vert\textbf{i} - \textbf{j}\Vert \leqslant  b_{t,n}}}^{}\sum_{t=1}^{T}\text{Cov}(p(A_{t,\textbf{i}}),p(A_{t,\textbf{j}})) + \sum_{\substack{\textbf{i},\textbf{j} \in \Lambda_{\mathbf{n}} \\ \Vert\textbf{i} - \textbf{j}\Vert > b_{t,n}}}^{}\sum_{t=1}^{T}\text{Cov}(p(A_{t,\textbf{i}}),p(A_{t,\textbf{j}}))\\
		&=&\mathcal{J}_3 + \mathcal{J}_4.
	\end{eqnarray*}
	Similarly to $\mathcal{J}_1$, we get:
	\begin{eqnarray*}
		\mathcal{J}_3 &=& \sum_{\substack{\textbf{i},\textbf{j} \in \Lambda_{\mathbf{n}} \\ \Vert\textbf{i} - \textbf{j}\Vert \leqslant  b_{t,n}}}^{}\sum_{t=1}^{T}\text{Cov}(p(A_{t,\textbf{i}}),p(A_{t,\textbf{j}}))\\
		&=& \sum_{\textbf{i} \in \Lambda_{\mathbf{n}}}^{}\sum_{\textbf{j} \in \Lambda_{\mathbf{n}} : \Vert\textbf{i}-\textbf{j}\Vert \leqslant  b_{t,n}}^{}\sum_{t=1}^{T}\text{Cov}(p(A_{t,\textbf{i}}),p(A_{t,\textbf{j}}))
		=  O\left(\widehat{\textbf{n}}Tb^d_{t,n}\widetilde{h}^2_{T,\textbf{n}}\right).
	\end{eqnarray*}
	Similarly to $\mathcal{J}_2$, under hypothesis \hyperref[H2]{H2}, by \cite[Corollary 3.3]{hwang13} and defining $t_{\textbf{j}}=\Vert\textbf{i}-\textbf{j}\Vert$ for each fixed $\textbf{i} \in \Lambda_{\mathbf{n}}$, we obtain:
	\begin{eqnarray*}
		\mathcal{J}_4 &=& \sum_{\textbf{i} \in \Lambda_{\mathbf{n}}}^{}\sum_{\textbf{j} \in \Lambda_{\mathbf{n}} : \Vert\textbf{i}-\textbf{j}\Vert > b_{t,n}}^{}\sum_{t=1}^{T}\text{Cov}(p(A_{t,\textbf{i}}),p(A_{t,\textbf{j}}))\\
		&\leqslant & \sum_{\textbf{i} \in \Lambda_{\mathbf{n}}}^{}\sum_{\textbf{j} \in \Lambda_{\mathbf{n}} : \Vert\textbf{i}-\textbf{j}\Vert > b_{t,n}}^{}\sum_{t=1}^{T}\left|\text{Cov}(p(A_{t,\textbf{i}}),p(A_{t,\textbf{j}}))\right|\\
		&\leqslant & C\widehat{\textbf{n}}T\left\Vert p(A_{t,\textbf{i}})\right\Vert_4\left\Vert p(A_{t,\textbf{j}})\right\Vert_4\sum_{\textbf{j} \in \Lambda_{\mathbf{n}} : t_{\textbf{j}} > b_{t,n}}^{}(\eta_2(t_{\textbf{j}}))^{1/2}\\
		&\leqslant & C\widehat{\textbf{n}}T\left\Vert p(A_{t,\textbf{i}})\right\Vert_4\left\Vert p(A_{t,\textbf{j}})\right\Vert_4\sum_{\textbf{j} \in \Lambda_{\mathbf{n}} : t_{\textbf{j}} > b_{t,n}}^{}t^{-\theta/2}_{\textbf{j}}
		= O\left(\widehat{\textbf{n}} Tb^{1-\theta/2}_{t,n}\widetilde{h}^{1/2}_{T,\textbf{n}}\right).
	\end{eqnarray*}
	Setting $b_{t,n} = \left\lfloor(\widetilde{h}_{T,\textbf{n}})^{-v/d}\right\rfloor$ with $v > d/(\theta - 2)$, and for $\theta > d + 1$, we obtain: 
	$
	\mathcal{I}_3 = O\left(\widehat{\textbf{n}}T\widetilde{h}_{T,\textbf{n}}\right).
	$
	Analogously to all the above, we get:
	\begin{eqnarray*}
		\mathcal{I}_4 &=& \sum_{\substack{\textbf{i},\textbf{j} \in \Lambda_{\mathbf{n}} \\ \textbf{i} \ne \textbf{j}}}^{}\sum_{s\ne t}^{T}\text{Cov}(p(A_{t,\textbf{i}}),p(A_{s,\textbf{j}}))
		= \sum_{\substack{\textbf{i},\textbf{j} \in \Lambda_{\mathbf{n}} \\ \Vert\textbf{i} - \textbf{j}\Vert \leqslant  b_{t,n}}}^{}\sum_{s= 1}^{T}\sum^{T}_{\substack{t = 1\\ |s - t| \leqslant  a_{n,t}}}\text{Cov}(p(A_{t,\textbf{i}}),p(A_{s,\textbf{j}}))\\
		&+& \sum_{\substack{\textbf{i},\textbf{j} \in \Lambda_{\mathbf{n}} \\ \Vert\textbf{i} - \textbf{j}\Vert \leqslant  b_{t,n}}}^{}\sum_{s = 1}^{T}\sum^{T}_{\substack{t = 1\\ |s- t| > a_{n,t}}}\text{Cov}(p(A_{t,\textbf{i}}),p(A_{s,\textbf{j}}))
		+ \sum_{\substack{\textbf{i},\textbf{j} \in \Lambda_{\mathbf{n}} \\ \Vert\textbf{i} - \textbf{j}\Vert > b_{t,n}}}^{}\sum_{s= 1}^{T}\sum^{T}_{\substack{t = 1\\ |s - t| \leqslant  a_{n,t}}}\text{Cov}(p(A_{t,\textbf{i}}),p(A_{s,\textbf{j}}))\\
		&+& \sum_{\substack{\textbf{i},\textbf{j} \in \Lambda_{\mathbf{n}} \\ \Vert\textbf{i} - \textbf{j}\Vert > b_{t,n}}}^{}\sum_{s= 1}^{T}\sum^{T}_{\substack{t = 1\\ |s - t| > a_{n,t}}}\text{Cov}(p(A_{t,\textbf{i}}),p(A_{s,\textbf{j}}))\\
		&=& \mathcal{J}_5 + \mathcal{J}_6 + \mathcal{J}_7 + \mathcal{J}_8.
	\end{eqnarray*}
	Using the same arguments as for $\mathcal{J}_1$ and $\mathcal{J}_3$, we obtain
	$
	\mathcal{J}_5 = O\!\left(\widehat{\mathbf{n}}T\,b_{t,n}^{d}a_{t,n}\widetilde{h}_{T,\mathbf{n}}^{2}\right).
	$
	Similarly, by reasoning as in the treatment of $\mathcal{J}_2$ and $\mathcal{J}_3$, it follows that
	$
	\mathcal{J}_6 = O\!\left(\widehat{\mathbf{n}}T\,a_{t,n}^{1-\theta/2}b_{t,n}^{d}\widetilde{h}_{T,\mathbf{n}}^{1/2}\right).
	$
	In the same way, arguing as for $\mathcal{J}_1$ and $\mathcal{J}_4$, we derive
	$
	\mathcal{J}_7 = O\!\left(\widehat{\mathbf{n}}T\,b_{t,n}^{1-\theta/2}a_{t,n}\widetilde{h}_{T,\mathbf{n}}^{1/2}\right),
	$
	whereas an argument similar to that used for $\mathcal{J}_2$ and $\mathcal{J}_4$ yields
	$
	\mathcal{J}_8 = O\!\left(\widehat{\mathbf{n}}T\,b_{t,n}^{1-\theta/2}a_{t,n}^{1-\theta/2}\widetilde{h}_{T,\mathbf{n}}^{1/2}\right).
	$
	
	Now, choosing
	\[
	a_{t,n}=\left\lfloor \widetilde{h}_{T,\mathbf{n}}^{-w/2}\right\rfloor
	\quad\text{with}\quad
	w>\frac{1}{\theta-2},
	\]
	and
	\[
	b_{t,n}=\left\lfloor \widetilde{h}_{T,\mathbf{n}}^{-v/(2d)}\right\rfloor
	\quad\text{with}\quad
	v>\frac{d}{\theta-2},
	\]
	we obtain, for $\theta>d+1$,
	$
	\mathcal{I}_4 = O\!\left(\widehat{\mathbf{n}}T\widetilde{h}_{T,\mathbf{n}}\right).
	$
	Therefore, we conclude that
	$$
	Var\left(\displaystyle\sum_{\textbf{i} \in \Lambda_{\mathbf{n}}}^{}\sum_{t=1}^{T}p(A_{t,\textbf{i}})\right) = O\left(\widehat{\textbf{n}}T\widetilde{h}_{T,\textbf{n}}\right).
	$$
\end{proof}
	\noindent
	In the following, we consider the space--time process \( X_{t,\mathbf{i}} \), which is \(\tau\)-weakly dependent, as defined in Section~\ref{sec3}.
	
	\begin{lemma}[]\label{lm2}
		Under assumptions \hyperref[H1]{H1}, \hyperref[H2]{H2}, \hyperref[H5]{H5}, and \hyperref[H6]{H6}, for all \( x \in \mathcal{Z} \), we have
		\[
		\left| \widetilde{H}_{T,\mathbf{n}}(x) - H(x) \right|
		= O\!\left( (\widetilde{h}_{T,\mathbf{n}})^2
		+ \sqrt{\dfrac{\log(\widehat{\mathbf{n}} T)}{\widehat{\mathbf{n}} T}} \right)
		\hspace{1cm}\text{a.s.}
		\]
	\end{lemma}
	
	\begin{proof}[Proof of Lemma~\ref{lm2}]
		Let us first consider \( \mathbb{E}\left[\widetilde{H}_{T,\textbf{n}}(x)\right] - H(x) \). We have:
		\begin{eqnarray*}
			\mathbb{E}\left[\widetilde{H}_{T,\textbf{n}}(x)\right] &=& \mathbb{E}\left[ \dfrac{1}{\widehat{\textbf{n}}T}\displaystyle\sum_{\textbf{i} \in \Lambda_{\mathbf{n}}}^{}\sum_{t=1}^{T}L\left(\dfrac{x-X_{t,\textbf{i}}}{\widetilde{h}_{T,\textbf{n}}}\right)\right] 
			= \displaystyle\int_{-\infty}^{+\infty}L\left(\dfrac{x-z}{\widetilde{h}_{T,\textbf{n}}}\right)h(z)\,dz.
		\end{eqnarray*}
		Integration by parts yields:
		\begin{eqnarray*}
			\mathbb{E}\left[\widetilde{H}_{T,\textbf{n}}(x)\right] &=& \left[L\left(\dfrac{x-z}{\widetilde{h}_{T,\textbf{n}}}\right)H(z)\right]_{-\infty}^{+\infty} + \dfrac{1}{\widetilde{h}_{T,\textbf{n}}}\displaystyle\int_{-\infty}^{+\infty}L'\left(\dfrac{x-z}{\widetilde{h}_{T,\textbf{n}}}\right)H(z)\,dz \\
			&=& 0 + \dfrac{1}{\widetilde{h}_{T,\textbf{n}}}\displaystyle\int_{-\infty}^{+\infty}L'\left(\dfrac{x-z}{\widetilde{h}_{T,\textbf{n}}}\right)H(z)\,dz.
		\end{eqnarray*}
		\noindent The term \( \left[L\left(\dfrac{x-z}{\widetilde{h}_{T,\textbf{n}}}\right)H(z)\right]_{-\infty}^{+\infty} \) vanishes because, at \( z = +\infty \), we have \( L\left(\dfrac{x-z}{\widetilde{h}_{T,\textbf{n}}}\right) = 0 \) and \( H(z) = 1 \), while at \( z = -\infty \), we have \( H(z) = 0 \) and \( L\left(\dfrac{x-z}{\widetilde{h}_{T,\textbf{n}}}\right) = 1 \). By setting \(\dfrac{x-z}{\widetilde{h}_{T,\textbf{n}}} = u\), we obtain:
		\[
		\mathbb{E}\left[\widetilde{H}_{T,\textbf{n}}(x)\right] = \displaystyle\int_{-\infty}^{+\infty}L'\left(u\right)H(x-u\widetilde{h}_{T,\textbf{n}})\, du.
		\]
		Applying the Taylor expansion with integral remainder, we obtain
		\[
		H\!\left(x - u\,\widetilde h_{T,\mathbf n}\right)
		= H(x) - u\,\widetilde h_{T,\mathbf n}\,H'(x)
		+ \big(u\,\widetilde h_{T,\mathbf n}\big)^{2}\!
		\int_{0}^{1} (1-\mu)\,H''\!\left(x - \mu\,u\,\widetilde h_{T,\mathbf n}\right)\,d\mu.
		\]
		
		\noindent Under assumptions \hyperref[H2]{H2} and \hyperref[H3]{H3}, since \(L'(u) = K(u)\) and the derivative \(H'' = h'\) is bounded, after simplification we obtain
		\begin{equation}
			\left|\,\mathbb{E}\!\left[\widetilde{H}_{T,\mathbf n}(x)\right] - H(x)\,\right|
			= O\!\left(\widetilde{h}_{T,\mathbf n}^{\,2}\right).
			\label{eq1}
		\end{equation}
		\noindent Moreover,
		\begin{eqnarray*}
			\widetilde{H}_{T,\textbf{n}}(x) - \mathbb{E}\left[\widetilde{H}_{T,\textbf{n}}(x)\right] &=&  \dfrac{1}{\widehat{\textbf{n}}T}\displaystyle\sum_{\textbf{i} \in \Lambda_{\mathbf{n}}}^{}\sum_{t=1}^{T}L\left(\dfrac{x-X_{t,\textbf{i}}}{\widetilde{h}_{T,\textbf{n}}}\right)  - \mathbb{E}\left[ \dfrac{1}{\widehat{\textbf{n}}T}\displaystyle\sum_{\textbf{i} \in \Lambda_{\mathbf{n}}}^{}\sum_{t=1}^{T}L\left(\dfrac{x-X_{t,\textbf{i}}}{\widetilde{h}_{T,\textbf{n}}}\right)\right] \\
			&=& \displaystyle\sum_{\textbf{i} \in \Lambda_{\mathbf{n}}}^{}\sum_{t=1}^{T} \Delta_{t,\textbf{i}} = G_{T,\textbf{n}}(x),
		\end{eqnarray*}
		where \( \Delta_{t,\textbf{i}} = \dfrac{1}{\widehat{\textbf{n}}T}\left(L\left(\dfrac{x-X_{t,\textbf{i}}}{\widetilde{h}_{T,\textbf{n}}}\right) - \mathbb{E}\left[L\left(\dfrac{x-X_{t,\textbf{i}}}{\widetilde{h}_{T,\textbf{n}}}\right)\right]\right) \). We represent the random variables \( \Delta_{t,\textbf{i}} \) in blocks as in \cite{amiri18} and \cite{wang231}. Considering the sequences \(r_k = \dfrac{n_k}{p_n + q_n}\), \(k = 1,\dots,d\), and \(r_{d+1} = \dfrac{T}{p_T + q_T}\), we define the block sums as follows:
		\begin{align*}
			U(1,n,\textbf{j},T,u,x) &= \sum_{\substack{\textbf{i} : i_k = j_k(p_n+q_n) + 1 \\ k = 1,\dots,d}}^{j_k(p_n + q_n) + p_n} \sum_{t = u(p_T + q_T) + 1}^{u(p_T + q_T) + p_T} \Delta_{t,\textbf{i}} \\
			U(2,n,\textbf{j},T,u,x) &= \sum_{\substack{\textbf{i} : i_k = j_k(p_n+q_n) + 1 \\ k = 1,\dots,d}}^{j_k(p_n + q_n) + p_n} \sum_{t = u(p_T + q_T) + p_T + 1}^{(u+1)(p_T + q_T)} \Delta_{t,\textbf{i}} \\
			U(3,n,\textbf{j},T,u,x) &= \sum_{\substack{\textbf{i} : i_k = j_k(p_n+q_n) + 1 \\ k = 1,\dots,d-1}}^{j_k(p_n + q_n) + p_n} \sum_{\substack{i_d = j_d(p_n+q_n) + p_n + 1}}^{(j_d + 1)(p_n + q_n)} \sum_{t = u(p_T + q_T) + 1}^{u(p_T + q_T) + p_T} \Delta_{t,\textbf{i}}
		\end{align*}
		and so on. We also have:
		\begin{align*}
			U(2^d,n,\textbf{j},T,u,x) &= \sum_{\substack{\textbf{i} : i_k = j_k(p_n+q_n) + p_n + 1 \\ k = 1,\dots,d}}^{j_k(p_n + q_n) + p_n} \sum_{t = u(p_T + q_T) + 1}^{u(p_T + q_T) + p_T} \Delta_{t,\textbf{i}} \\
			U(2^{d+1},n,\textbf{j},T,u,x) &= \sum_{\substack{\textbf{i} : i_k= j_k(p_n+q_n) + p_n + 1 \\ k = 1,\dots,d}}^{(j_k + 1)(p_n + q_n)} \sum_{t = u(p_T + q_T) + p_T + 1}^{(u+1)(p_T + q_T)} \Delta_{t,\textbf{i}}.
		\end{align*}
		For any integer \(1  \leqslant  m  \leqslant  2^d\), we define:
		\begin{align*}
			E(n,m,T,x) = \sum_{\substack{j_k = 0 \\ k = 1,\dots,d}}^{r_k - 1} \sum_{u = 0}^{r_{d+1} - 1} U(m,n,\textbf{j},T,u,x).
		\end{align*}
		We may then write \[G_{T,\textbf{n}}(x)= \displaystyle\sum_{m = 1}^{2^{d+1}}E(n,m,T,x)\].
		\noindent Note that for all \(\epsilon > 0\),
		\[
		\begin{aligned}
			\mathbb{P}\left(\left|G_{T,\textbf{n}}(x)\right| > \epsilon\right)
			&= \mathbb{P}\left(
			\left|\sum_{m = 1}^{2^{d+1}} E(n,m,T,x)\right| > \epsilon
			\right) \\
			&\leqslant 2^{d+1}\,
			\mathbb{P}\left(
			\left|E(n,m,T,x)\right| > \frac{\epsilon}{2^{d+1}}
			\right).
		\end{aligned}
		\]
		Without loss of generality, for a representative block \(E(n,1,T,x)\), we have:
		\begin{eqnarray*}
			U(1,n,\textbf{j},T,u,x) &\leqslant & \left\|U(1,n,\textbf{j},T,u,x)\right\|_\infty \\
			&=& \left\|\sum_{\substack{\textbf{i} : i_k = j_k(p_n+q_n) + 1 \\ k = 1,...,d}}^{j_k(p_n + q_n) + p_n} \sum_{t=u(p_T + q_T) + 1}^{u(p_T + q_T) + p_T} \Delta_{t,\textbf{i}}\right\|_{\infty} \\
			&\leqslant & \sum_{\substack{\textbf{i} : i_k = j_k(p_n+q_n) + 1 \\ k = 1,...,d}}^{j_k(p_n + q_n) + p_n} \sum_{t=u(p_T + q_T) + 1}^{u(p_T + q_T) + p_T} \left\|\Delta_{t,\textbf{i}}\right\|_\infty \\
			&\leqslant & p_n^d p_T \left\|\Delta_{t,\textbf{i}}\right\|_\infty \\
			&\leqslant & \frac{p_n^d p_T \left\|L\right\|_{\infty}}{\widehat{\textbf{n}}T} =: b_{T,\textbf{n}}.
		\end{eqnarray*}
		Following \cite{amiri18}, we enumerate the variables \(U(1,n,\textbf{j},T,u,x)\) in an arbitrary order and denote them by \(\hat{U}_1, ..., \hat{U}_M\), where \(M = r_1 \cdots r_{d+1}\). Then, using \cite[Lemma 5]{dedecker04}, we couple the \(\tau\)-weakly dependent variables \(\hat{U}_i\) with independent counterparts \(\tilde{U}_i\), and apply Markov's inequality followed by Bernstein's inequality. For a sequence \(\epsilon_{T,\mathbf{n}} \to 0\), we obtain:
		\begin{equation*}
			\begin{aligned}
				\mathbb{P}\left( \left|E(n,1,T,x)\right| > \epsilon_{T,\mathbf{n}} \right)
				&= \mathbb{P}\left(
				\left| \hat{H}_{T,\mathbf{n}}(x)
				- \mathbb{E}\left[\hat{H}_{T,\mathbf{n}}(x)\right] \right|
				> \epsilon_{T,\mathbf{n}}
				\right) \\
				&\leqslant \frac{2M}{\epsilon_{T,\mathbf{n}}}
				\, b_{T,\mathbf{n}} \tau_1(q_n)\tau_2(q_T) \\
				&\quad + 2 \exp\left(
				-\frac{\epsilon_{T,\mathbf{n}}^2}{
					8\,\operatorname{Var}\left(\widetilde{H}_{T,\mathbf{n}}(x)\right)
					+ \frac{4}{3} b_{T,\mathbf{n}} \epsilon_{T,\mathbf{n}}
				}
				\right) \\
				&=: \mathcal{A}_1 + \mathcal{A}_2.
			\end{aligned}
		\end{equation*}
		Here, \(q_n\) denotes the minimal spatial distance between two blocks, and \(q_T\) the minimal temporal distance. Choosing
		\[
		p_n = \left\lfloor\left( \frac{\hat{\mathbf{n}}}{\log \hat{\mathbf{n}}} \right)^{1/2d}\right\rfloor, \quad
		p_T = \left\lfloor\left( \frac{T}{\log T} \right)^{1/2}\right\rfloor, \quad
		\epsilon_{T,\mathbf{n}} = \sqrt{ \frac{\log(\hat{\mathbf{n}}T)}{ \hat{\mathbf{n}}T } },
		\]
		and using the variance bounds from Lemma~\ref{lm1}, where \(\operatorname{Var}\left( \widetilde{H}_{T,\textbf{n}} \right) = O\left( \frac{ \widetilde{h}_{T,\textbf{n}} }{ \widehat{\textbf{n}}T } \right)\), we have:
		
		\begin{eqnarray*}
			\mathcal{A}_2 &=& 2 \exp\left( -\frac{ \epsilon^2_{T,\textbf{n}} }{ 8 \operatorname{Var}\left( \widetilde{H}_{T,\textbf{n}}(x) \right) + \frac{4}{3} b_{T,\textbf{n}} \epsilon_{T,\textbf{n}} } \right) \\
			&\leqslant & 2 \exp\left( -\frac{ \log(\widehat{\textbf{n}}T) }{ 8C \widetilde{h}_{T,\textbf{n}} + \frac{4 \left\|L\right\|_{\infty} (\log(\widehat{\textbf{n}}T))^{1/2} }{ 3 (\log \widehat{\textbf{n}} \log T)^{1/2} } } \right) \\
			&\leqslant & 2 \exp\left( -B_{T,\textbf{n}} \log(\widehat{\textbf{n}}T) \right)
			= \frac{2}{(T \widehat{\textbf{n}})^{B_{T,\textbf{n}}}} := \alpha_{T,\textbf{n}},
		\end{eqnarray*}
		where the series with general term \(\alpha_{T,\textbf{n}}\) converges since \(B_{T,\textbf{n}} \to \infty\).
		For \(\mathcal{A}_1\), noting that \(M \leqslant  \frac{T \widehat{\textbf{n}}}{p_n^d p_T}\), we obtain:
		\begin{eqnarray*}
			\mathcal{A}_1
			&=& \frac{2M}{\epsilon_{T,\textbf{n}}} b_{T,\textbf{n}} \tau_1(q_n) \tau_2(q_T)
			\leqslant  \frac{2}{\epsilon_{T,\textbf{n}}} \tau_1(q_n) \tau_2(q_T)
			\leqslant  \frac{2(\widehat{\textbf{n}}T)^{1/2}}{(\log \widehat{\textbf{n}}T)^{1/2}} \tau_1(q_n) \tau_2(q_T).
		\end{eqnarray*}
		We now choose \(q_n = \left\lfloor\log(\widehat{\textbf{n}})\right\rfloor\) and \(q_T = \left\lfloor\log(T)\right\rfloor\). According to assumption \hyperref[H1]{H1}, we have:
		\[
		\tau_1(q_n) = O\left( \frac{1}{\widehat{\textbf{n}}^{\rho_1}} \right), \quad \tau_2(q_T) = O\left( \frac{1}{T^{\rho_2}} \right),
		\]
		so that:
		\[
		\mathcal{A}_1 \leqslant  \frac{2C}{ (\log \widehat{\textbf{n}}T)^{1/2} T^{\rho_2 - 1/2} \widehat{\textbf{n}}^{\rho_1 - 1/2} } =: \beta_{T,\textbf{n}}.
		\]
		It is enough to choose \(\rho_1 > \frac{3}{2}\) and \(\rho_2 > \frac{3}{2}\) to ensure that the series with general term \(\beta_{T,\textbf{n}}\) converges. Then, by the Borel--Cantelli lemma, we conclude that
		\begin{eqnarray}
			\left| \widetilde{H}_{T,\textbf{n}}(x) - \mathbb{E}\left[ \widetilde{H}_{T,\textbf{n}}(x) \right] \right| = O\left( \sqrt{ \frac{ \log(\widehat{\textbf{n}}T) }{ \widehat{\textbf{n}}T } } \right) \quad \text{a.s.}
			\label{eq2}
		\end{eqnarray}
		Combining \eqref{eq1} and \eqref{eq2} yields
		\[
		\left| \widetilde{H}_{T,\textbf{n}}(x) - H(x) \right| = O\left( \widetilde{h}_{T,\textbf{n}}^2 + \sqrt{ \frac{ \log(\widehat{\textbf{n}}T) }{ \widehat{\textbf{n}}T } } \right) \quad \text{a.s.}
		\]
	\end{proof}
	
	\begin{lemma}\label{lm3}
		Under assumptions \hyperref[H1]{H1}, \hyperref[H2]{H2}, \hyperref[H5]{H5}, and \hyperref[H6]{H6}, we have:
		\begin{eqnarray*}
			\displaystyle\sup_{x \in \mathcal{C}}\left|\widetilde{H}_{T,\textbf{n}}(x) - H(x)\right| = O\left((\widetilde{h}_{T,\textbf{n}})^2 + \sqrt{\dfrac{\log(\widehat{\textbf{n}}T)}{\widehat{\textbf{n}}T}}\right) \hspace{1cm}\text{a.s.}
		\end{eqnarray*}
		\label{l3}
	\end{lemma}
	
	\begin{proof}[Proof of Lemma~\ref{lm3}]
		\noindent We begin by examining the expression
		\[
		\sup_{x \in \mathcal{C}} \left| \mathbb{E}\!\left[\widetilde{H}_{T,\textbf{n}}(x)\right] - H(x) \right|.
		\]
		Based on the proof of Lemma~\ref{lm2}, we obtain that, for all \( x \in \mathcal{Z} \), there exists a constant \( C > 0 \) such that
		\[
		\left| \mathbb{E}\!\left[\widetilde{H}_{T,\textbf{n}}(x)\right] - H(x) \right| \leqslant C \, (\widetilde{h}_{T,\textbf{n}})^2.
		\]
		Since \( \mathcal{C} \) is a compact subset of \( \mathcal{Z} \), it follows immediately that
		
		\begin{equation}
			\sup_{x \in \mathcal{C}} \left| \mathbb{E}\left[\widetilde{H}_{T,\textbf{n}}(x)\right] - H(x) \right| = O\left( \widetilde{h}_{T,\textbf{n}}^2 \right).
			\label{eq3}
		\end{equation}
		We now turn to the study of
		\[
		\sup_{x \in \mathcal{C}} \left| \widetilde{H}_{T,\textbf{n}}(x) - \mathbb{E}\left[ \widetilde{H}_{T,\textbf{n}}(x) \right] \right|.
		\]
		Since \( \mathcal{C} \subset \mathcal{Z} \subset \mathbb{R} \) is compact, it can be covered by a finite number \( v \) of intervals \( I_k \), each centered at some point \( x_k \) and of radius \( r \), such that for every \( x \in \mathcal{C} \), there exists \( k \in \{1, \dots, v\} \) with \( x \in I_k \).
		Thus,
		\[
		\sup_{x \in \mathcal{C}} \left| \widetilde{H}_{T,\textbf{n}}(x) - \mathbb{E}\left[ \widetilde{H}_{T,\textbf{n}}(x) \right] \right|
		= \max_{1  \leqslant  k  \leqslant  v} \sup_{x \in I_k} \left| \widetilde{H}_{T,\textbf{n}}(x) - \mathbb{E}\left[ \widetilde{H}_{T,\textbf{n}}(x) \right] \right|.
		\]
		
		For every \( x \in I_k \), the triangle inequality gives:
		\begin{eqnarray*}
			\left| \widetilde{H}_{T,\textbf{n}}(x) - \mathbb{E}\left[\widetilde{H}_{T,\textbf{n}}(x)\right] \right|
			& \leqslant &
			\left| \widetilde{H}_{T,\textbf{n}}(x) - \widetilde{H}_{T,\textbf{n}}(x_k) \right| \\
			&& + \left| \widetilde{H}_{T,\textbf{n}}(x_k) - \mathbb{E}\left[ \widetilde{H}_{T,\textbf{n}}(x_k) \right] \right| \\
			&& + \left| \mathbb{E}\left[ \widetilde{H}_{T,\textbf{n}}(x_k) \right] - \mathbb{E}\left[ \widetilde{H}_{T,\textbf{n}}(x) \right] \right|.
		\end{eqnarray*}
		Taking the supremum over each \( I_k \) and then the maximum over \( k \), we obtain:
		\begin{equation}
			\sup_{x \in \mathcal{C}} \left| \widetilde{H}_{T,\textbf{n}}(x) - \mathbb{E}\left[ \widetilde{H}_{T,\textbf{n}}(x) \right] \right|  \leqslant  Q_1 + Q_2 + Q_3,
			\label{eq4}
		\end{equation}
		where
		\[
		Q_1 = \max_{1  \leqslant  k  \leqslant  v} \sup_{x \in I_k} \left| \widetilde{H}_{T,\textbf{n}}(x) - \widetilde{H}_{T,\textbf{n}}(x_k) \right|, \quad
		Q_2 = \max_{1  \leqslant  k  \leqslant  v} \left| \widetilde{H}_{T,\textbf{n}}(x_k) - \mathbb{E}\left[ \widetilde{H}_{T,\textbf{n}}(x_k) \right] \right|,
		\]
		\[
		Q_3 = \max_{1  \leqslant  k  \leqslant  v} \sup_{x \in I_k} \left| \mathbb{E}\left[ \widetilde{H}_{T,\textbf{n}}(x) \right] - \mathbb{E}\left[ \widetilde{H}_{T,\textbf{n}}(x_k) \right] \right|.
		\]
		We now treat \( Q_1 \), \( Q_2 \), and \( Q_3 \) individually. Thanks to the Lipschitz continuity of the kernel \( L \), with Lipschitz constant \( \mu \), we have for any \( x \in I_k \),
		\[
		\left| L\left( \frac{x - X_{t,\textbf{i}} }{\widetilde{h}_{T,\textbf{n}}} \right) - L\left( \frac{x_k - X_{t,\textbf{i}}}{\widetilde{h}_{T,\textbf{n}}} \right) \right|  \leqslant  \frac{ \mu |x - x_k| }{ \widetilde{h}_{T,\textbf{n}} }.
		\]
		Since \( |x - x_k|  \leqslant  r \) on \( I_k \), choosing
		$
		r = \widetilde{h}_{T,\textbf{n}} \sqrt{ \frac{ \log(\widehat{\textbf{n}}T) }{ \widehat{\textbf{n}}T } },
		$
		we get
		$
		Q_1 = O\left( \sqrt{ \frac{ \log(\widehat{\textbf{n}}T) }{ \widehat{\textbf{n}}T } } \right).
		$
		By an argument analogous to the one used for \( Q_1 \), we also obtain
		$
		Q_3 = O\left( \sqrt{ \frac{ \log(\widehat{\textbf{n}}T) }{ \widehat{\textbf{n}}T } } \right).
		$
		For controlling \( Q_2 \), for any \( \epsilon > 0 \),
		$
		\mathbb{P}\left( Q_2 > \epsilon \right)  \leqslant  v \sup_{x \in \mathcal{C}} \mathbb{P}\left( \left| \widetilde{H}_{T,\textbf{n}}(x) - \mathbb{E}\left[ \widetilde{H}_{T,\textbf{n}}(x) \right] \right| > \epsilon \right).
		$
		However, we have
		
		\begin{equation}
			\left| \widetilde{H}_{T,\textbf{n}}(x) - \mathbb{E}\left[ \widetilde{H}_{T,\textbf{n}}(x) \right] \right| = O\left( \sqrt{ \frac{ \log(\widehat{\textbf{n}}T) }{ \widehat{\textbf{n}}T } } \right) \quad \text{a.s.}
			\label{eq5}
		\end{equation}
		Combining results \eqref{eq3}, \eqref{eq4}, and \eqref{eq5}, we finally obtain:
		\[
		\sup_{x \in \mathcal{C}} \left| \widetilde{H}_{T,\textbf{n}}(x) - H(x) \right| = O\left( \widetilde{h}_{T,\textbf{n}}^2 + \sqrt{ \frac{ \log(\widehat{\textbf{n}}T) }{ \widehat{\textbf{n}}T } } \right) \quad \text{a.s.}
		\]
	\end{proof}
	
	\begin{lemma}\label{lm4}
		Under assumptions \hyperref[H1]{H1}, \hyperref[H2]{H2}, \hyperref[H5]{H5}, and \hyperref[H6]{H6}, for every $x \in \mathcal{Z}$, we have:
		\begin{eqnarray*}
			\left|\widetilde{f}^{*}_{T,\textbf{n}}(x) - f^*(x)\right| = O\left((\widetilde{h}_{T,\textbf{n}})^{2}+\sqrt{\dfrac{\log( \widehat{\textbf{n}}T)}{\widehat{\textbf{n}}T\widetilde{h}_{T,\textbf{n}}}}\right) \hspace{1cm} \text{p.s.}
		\end{eqnarray*}
		\label{l4}
	\end{lemma}
	
	\begin{proof}[Proof of Lemma \ref{lm4}]
		\noindent Consider the following estimator:
		\[
		\widetilde{f}^{*}_{T,\textbf{n}}(x) = \frac{1}{\widehat{\textbf{n}}T \widetilde{h}_{T,\textbf{n}}} \sum_{\textbf{i} \in \Lambda_{\mathbf{n}}} \sum_{t=1}^{T} \delta_{t,\textbf{i}} K\left( \frac{x - X_{t,\textbf{i}}}{\widetilde{h}_{T,\textbf{n}}} \right).
		\]
		We want to show that \(\widetilde{f}^{*}_{T,\textbf{n}}(x)\) converges to \(f^*(x)\) for all \(x \in \mathcal{Z}\).
		Let us begin by analyzing the bias term \(\mathbb{E}[\widetilde{f}^{*}_{T,\textbf{n}}(x)] - f^*(x)\).
		We have:
		
		\begin{eqnarray*}
			\mathbb E\!\left[\delta_{t,\mathbf{i}}\,K\!\left(\frac{x-X_{t,\mathbf{i}}}{\widetilde h_{T,\mathbf n}}\right)\right]
			&=& \mathbb E\!\left[\; \mathbb E\!\left[\delta_{t,\mathbf{i}}\,K\!\left(\frac{x-X_{t,\mathbf{i}}}{\widetilde h_{T,\mathbf n}}\right)\,\Bigm|\,\delta_{t,\mathbf{i}}\right]\;\right] \\[1mm]
			&=& \sum_{d\in\{0,1\}}\mathbb P(\delta_{t,\mathbf{i}}=d)\;
			\underbrace{\mathbb E\!\left[\delta_{t,\mathbf{i}}\,K\!\left(\frac{x-X_{t,\mathbf{i}}}{\widetilde h_{T,\mathbf n}}\right)\,\Bigm|\,\delta_{t,\mathbf{i}}=d\right]}_{\text{computed for }d=0,1}\\[1mm]
			&=& \sum_{d\in\{0,1\}}\mathbb P(\delta_{t,\textbf{i}}=d)\;
			\underbrace{d\,\mathbb E\!\left[\,K\!\left(\frac{x-X_{t,\textbf{i}}}{\widetilde{h}_{T,\textbf{n}}}\right)\,\Bigm|\,\delta_{t,\textbf{i}}=d\right]}_{\delta_{t,\textbf{i}}=d\text{ is constant}} \\[1mm]
			&=& \mathbb P(\delta_{t,\textbf{i}}=1)\,\mathbb E\!\left[\,K\!\left(\frac{x-X_{t,\textbf{i}}}{\widetilde{h}_{T,\textbf{n}}}\right)\,\Bigm|\,\delta_{t,\textbf{i}}=1\right]\\
			\;&+&\; \underbrace{\mathbb P(\delta_{t,\textbf{i}}=0)\,\cdot 0}_{\text{the term for }\delta_{t,\textbf{i}}=0\text{ vanishes}} \\[1mm]
			&=& \mathbb P(\delta_{t,\textbf{i}}=1)\,\mathbb E\!\left[\,K\!\left(\frac{x-\mathcal{T}_{t,\mathbf{i}}}{\widetilde{h}_{T,\textbf{n}}}\right)\,\Bigm|\,\delta_{t,\textbf{i}}=1\right]\\
			&&\qquad(\text{since }\delta_{t,\textbf{i}}=1\Rightarrow X_{t,\textbf{i}}=\mathcal{T}_{t,\mathbf{i}}).
		\end{eqnarray*}
		\begin{align*}
			&\mathbb P(\delta_{t,\mathbf{i}}=1)\,
			\mathbb E\!\left[
			K\!\left(\frac{x-\mathcal{T}_{t,\mathbf{i}}}{\widetilde{h}_{T,\mathbf{n}}}\right)
			\,\Bigm|\,\delta_{t,\mathbf{i}}=1
			\right] \\
			&\quad=
			\mathbb P(\delta_{t,\mathbf{i}}=1)
			\int_{\mathbb R}
			K\!\left(\frac{x-t}{\widetilde{h}_{T,\mathbf{n}}}\right)\,
			f_{T\,|\,\delta_{t,\mathbf{i}}=1}(t)\,dt \\[1mm]
			&\quad=
			\mathbb P(\delta_{t,\mathbf{i}}=1)
			\int_{\mathbb R}
			K\!\left(\frac{x-t}{\widetilde{h}_{T,\mathbf{n}}}\right)\,
			\frac{
				f(t)\,\mathbb P(\delta_{t,\mathbf{i}}=1\,|\,\mathcal{T}_{t,\mathbf{i}}=t)
			}{
				\mathbb P(\delta_{t,\mathbf{i}}=1)
			}\,dt \\[1mm]
			&\quad=
			\int_{\mathbb R}
			K\!\left(\frac{x-t}{\widetilde{h}_{T,\mathbf{n}}}\right)\,
			f(t)\,
			\mathbb P(\delta_{t,\mathbf{i}}=1\,|\,\mathcal{T}_{t,\mathbf{i}}=t)\,dt.
		\end{align*}
	
	\begin{align*}
		&\quad=
		\int_{\mathbb R}
		K\!\left(\frac{x-t}{\widetilde{h}_{T,\mathbf{n}}}\right)\,
		f(t)\,
		\mathbb P(\mathcal{T}_{t,\mathbf{i}}\le C_{t,\mathbf{i}}\,|\,\mathcal{T}_{t,\mathbf{i}}=t)\,dt \\[1mm]
		&\quad=
		\int_{\mathbb R}
		K\!\left(\frac{x-t}{\widetilde{h}_{T,\mathbf{n}}}\right)\,
		f(t)\,
		\mathbb P(C_{t,\mathbf{i}}\ge t)\,dt \\[1mm]
		&\quad=
		\int_{\mathbb R}
		K\!\left(\frac{x-t}{\widetilde{h}_{T,\mathbf{n}}}\right)\,
		f(t)\,\bigl(1-G(t)\bigr)\,dt.
	\end{align*}
	
		With $u=\dfrac{x-t}{\widetilde{h}_{T,\mathbf{n}}}$, we have
		\[
		\mathbb E\!\left[\delta_{t,\mathbf{i}}\,K\!\left(\frac{x-X_{t,\mathbf{i}}}{\widetilde{h}_{T,\mathbf{n}}}\right)\right]
		= \frac{1}{\widetilde h_{T,\mathbf{n}}}\int_{\mathbb R} K(u)\,f^{*}\!\left(x-\widetilde{h}_{T,\mathbf{n}}\,u\right)\,du,
		\]
		\[
		\mathbb{E}\!\left[\widetilde f^{*}_{T,\mathbf n}(x)\right]
		= \int_{\mathbb{R}} K(u)\, f^*(x - \widetilde h_{T,\mathbf n}\,u)\,du.
		\]
		
		By a Taylor expansion with integral remainder, we have:
		\[
		f^*(x - \widetilde h_{T,\mathbf n}u)
		= f^*(x) - \widetilde h_{T,\mathbf n}u\,(f^*)'(x)
		+ (\widetilde h_{T,\mathbf n}u)^2
		\int_{0}^{1} (1-\mu)\,(f^*)''\!\big(x - \mu\,\widetilde h_{T,\mathbf n}u\big)\,d\mu.
		\]
		
		Thus, under assumption \hyperref[H2]{H2}, we obtain:
		\[
		\begin{aligned}
			\mathbb{E}\!\left[\widetilde f^{*}_{T,\mathbf n}(x)\right] - f^*(x)
			&= \int_D K(u)\Big(f^*(x - \widetilde h_{T,\mathbf n}u) - f^*(x)\Big)\,du \\[2pt]
			&= -\,\widetilde h_{T,\mathbf n}(f^*)'(x)\int_D u\,K(u)\,du\\
			&+ \int_D K(u)\,(\widetilde h_{T,\mathbf n}u)^2
			\int_0^1 (1-\mu)\,(f^*)''\!\big(x - \mu\,\widetilde h_{T,\mathbf n}u\big)\,d\mu\,du \\[2pt]
			&\le (\widetilde h_{T,\mathbf n})^2\int_D |u|^2K(u)\,
			\int_0^1 \left|(1-\mu)\,(f^*)''\!\big(x - \mu\,\widetilde h_{T,\mathbf n}u\big)\right|\,d\mu\,du.
		\end{aligned}
		\]
		Under \hyperref[H2]{H2} and \hyperref[H3]{H3}, it follows that:
		\begin{eqnarray}
			\label{eq6}
			\Big|\mathbb{E}\!\left[\widetilde f^{*}_{T,\mathbf n}(x)\right] - f^*(x)\Big|
			= O\!\big(\widetilde h_{T,\mathbf n}^2\big).
		\end{eqnarray}
		
		Let us now examine the deviation term \(\left| \widetilde{f}^{*}_{T,\textbf{n}}(x) - \mathbb{E}[ \widetilde{f}^{*}_{T,\textbf{n}}(x) ] \right|\).
		Following the same approach as in the proof of Lemma~\ref{lm2} replacing \(\widetilde{H}_{T,\textbf{n}}\) by \(\widetilde{f}^{*}_{T,\textbf{n}}\) and applying Lemma~\ref{lm1}, we obtain:
		\[
		\operatorname{Var} \left( \widetilde{f}^{*}_{T,\textbf{n}}(x) \right) = O\!\left( \frac{1}{\widehat{\textbf{n}}\, T \,\widetilde{h}_{T,\textbf{n}}} \right).
		\]
		and by choosing
		\[
		p_n = \left\lfloor\left( \frac{\widehat{\textbf{n}} (\widetilde{h}_{T,\textbf{n}})^{1/2} }{ \log \widehat{\textbf{n}} } \right)^{1/(2d)}\right\rfloor, \quad
		p_T = \left\lfloor\left( \frac{T (\widetilde{h}_{T,\textbf{n}})^{1/2} }{ \log T } \right)^{1/2}\right\rfloor, \quad
		\epsilon_{T,\textbf{n}} = \sqrt{ \frac{ \log(\widehat{\textbf{n}}T) }{ \widehat{\textbf{n}} T \widetilde{h}_{T,\textbf{n}} } },
		\]
		we obtain:
		\begin{equation}
			\left| \widetilde{f}^{*}_{T,\textbf{n}}(x) - \mathbb{E}[ \widetilde{f}^{*}_{T,\textbf{n}}(x) ] \right| = O\left( \sqrt{ \frac{ \log(\widehat{\textbf{n}}T) }{ \widehat{\textbf{n}} T \widetilde{h}_{T,\textbf{n}} } } \right) \quad \text{a.s.}
			\label{eq7}
		\end{equation}
		
		Combining \eqref{eq6} and \eqref{eq7}, we conclude that:
		\[
		\left| \widetilde{f}^{*}_{T,\textbf{n}}(x) - f^*(x) \right| = O\left( \widetilde{h}_{T,\textbf{n}}^2 + \sqrt{ \frac{ \log(\widehat{\textbf{n}}T) }{ \widehat{\textbf{n}} T \widetilde{h}_{T,\textbf{n}} } } \right) \quad \text{a.s.}
		\]
	\end{proof}
	
	\begin{lemma}\label{lm5}
		Under assumptions \hyperref[H1]{H1}, \hyperref[H2]{H2}, \hyperref[H5]{H5} and \hyperref[H6]{H6}, we have:
		\begin{eqnarray*}
			\sup_{x \in \mathcal{C}}\left|\widetilde{f}^{*}_{T,\textbf{n}}(x) - f^*(x)\right| = O\left(\widetilde{h}_{T,\textbf{n}}^{2}+\sqrt{\dfrac{\log( \widehat{\textbf{n}}T)}{\widehat{\textbf{n}}T\widetilde{h}_{T,\textbf{n}}}}\right) \hspace{1cm} \text{a.s.}
		\end{eqnarray*}
		\label{l5}
	\end{lemma}
	
	\begin{proof}[Proof of Lemma \ref{lm5}]
		Using results \eqref{eq6} and \eqref{eq7}, together with the reasoning followed in the proof of Lemma~\ref{l3}, we obtain:
		\[
		\sup_{x \in \mathcal{C}} \left| \widetilde{f}^*_{T,\textbf{n}}(x) - f^*(x) \right| = O\left( \widetilde{h}_{T,\textbf{n}}^2 + \sqrt{ \frac{ \log(\widehat{\textbf{n}}T) }{ \widehat{\textbf{n}}T \widetilde{h}_{T,\textbf{n}}} } \right) \quad \text{a.s.}
		\]
	\end{proof}

\begin{lemma}\label{lm6}
		Under assumptions \hyperref[H1]{H1}, \hyperref[H2]{H2}, \hyperref[H4]{H4} and \hyperref[H6]{H6}, we have:
		\begin{eqnarray*}
			\left|\widehat{H}^{}_{T,\textbf{n}}(x) - H(x)\right| = O\left(\left(\dfrac{k_{T,\textbf{n}}}{\widehat{\textbf{n}}T}\right)^2 +\sqrt{\dfrac{\log(\widehat{\textbf{n}}T)}{\widehat{\textbf{n}}T}}\right) \hspace{1cm}\text{a.s.}
		\end{eqnarray*}
	\end{lemma}
	
	\begin{proof}[Proof of Lemma~\ref{lm6}]
		Let $\beta \in ]0,1[$.
		Suppose that $D^-_{T,\mathbf{n}}(\beta) \leqslant h_{T,\mathbf{n}} \leqslant D^+_{T,\mathbf{n}}(\beta)$.
		Under assumption \hyperref[H2]{H2}, for all $u \in \mathbb{R}$ and all $t \in ]0,1[$, we have:
		
		\begin{itemize}
			\item If $u \geqslant 0$, then $tu < u$, hence $L(u) \geqslant L(tu)$.
			\item If $u < 0$, then $u < tu$, hence $L(u) \leqslant L(tu)$.
		\end{itemize}
		Hence, for $u \geqslant 0$:
		\begin{align*}
			L\left( \frac{x-X_{t,\mathbf{i}} }{D^+_{T,\mathbf{n}}(\beta)} \right)
			&\leqslant L\left( \frac{x-X_{t,\mathbf{i}}}{h_{T,\mathbf{n}}} \right)
			\leqslant L\left( \frac{x-X_{t,\mathbf{i}} }{D^-_{T,\mathbf{n}}(\beta)} \right),
		\end{align*}
		which implies:
		\[
		\widehat{H}^{+}_{T,\mathbf{n}}(x) \leqslant \widehat{H}_{T,\mathbf{n}}(x) \leqslant \widehat{H}^{-}_{T,\mathbf{n}}(x).
		\]
		For $u < 0$, the inequality is reversed:
		\begin{align*}
			L\left( \frac{x-X_{t,\mathbf{i}} }{D^-_{T,\mathbf{n}}(\beta)} \right)
			&\leqslant L\left( \frac{x-X_{t,\mathbf{i}}}{h_{T,\mathbf{n}}} \right)
			\leqslant L\left( \frac{x-X_{t,\mathbf{i}}}{D^+_{T,\mathbf{n}}(\beta)} \right),
		\end{align*}
		which leads to:
		\[
		\widehat{H}^{-}_{T,\mathbf{n}}(x) \leqslant \widehat{H}_{T,\mathbf{n}}(x) \leqslant \widehat{H}^{+}_{T,\mathbf{n}}(x),
		\]
		where:
		\begin{align*}
			\widehat{H}^{\pm}_{T,\mathbf{n}}(x) = \frac{1}{\widehat{\mathbf{n}} T} \sum_{\mathbf{i} \in \Lambda_{\mathbf{n}}} \sum_{t=1}^{T} L\left( \frac{x-X_{t,\mathbf{i}}}{D^{\pm}_{T,\mathbf{n}}(\beta)} \right).
		\end{align*}
		Thus, depending on the sign of $u$, the estimator $\widehat{H}_{T,\mathbf{n}}(x)$ is bounded between two kernel-based quantities controlled by $D^-_{T,\mathbf{n}}(\beta)$ and $D^+_{T,\mathbf{n}}(\beta)$.
		Let us now consider \(\widehat{H}^{-}_{T,\mathbf{n}}(x)\) (resp. \(\widehat{H}^{+}_{T,\mathbf{n}}(x)\)).
		Since \(D^-_{T,\mathbf{n}}(\beta)\) (resp. \(D^+_{T,\mathbf{n}}(\beta)\)) is deterministic, the proof of the almost sure convergence of \(\widehat{H}^{-}_{T,\mathbf{n}}(x)\) (resp. \(\widehat{H}^{+}_{T,\mathbf{n}}(x)\)) to \(H(x)\) follows exactly the same arguments as for \(\widetilde{H}_{T,\mathbf{n}}(x)\) in Lemma~\ref{lm2}.
		We choose:
		\[
		D^-_{T,\mathbf{n}}(\beta) = \frac{k_{T,\mathbf{n}}}{\hat{\mathbf{n}}T} \beta^{1/2},
		\quad
		D^+_{T,\mathbf{n}}(\beta) = \frac{k_{T,\mathbf{n}}}{\hat{\mathbf{n}}T} \beta^{-1/2},
		\]
		where \(k_{T,\mathbf{n}}\) is a sequence tending to infinity and satisfying assumption \hyperref[H4]{H4}.
		Then, we obtain:
		\[
		\left| \widehat{H}^{-}_{T,\mathbf{n}}(x) - H(x) \right| = O\left( \left( \frac{k_{T,\mathbf{n}}}{\hat{\mathbf{n}}T} \right)^2 + \sqrt{ \frac{ \log(\hat{\mathbf{n}}T) }{ \hat{\mathbf{n}}T } } \right) \quad \text{a.s.},
		\]
		and similarly,
		\[
		\left| \widehat{H}^{+}_{T,\mathbf{n}}(x) - H(x) \right| = O\left( \left( \frac{k_{T,\mathbf{n}}}{\hat{\mathbf{n}}T} \right)^2 + \sqrt{ \frac{ \log(\hat{\mathbf{n}}T) }{ \hat{\mathbf{n}}T } } \right) \quad \text{a.s.}
		\]
		It follows that:
		\[
		\left| \widehat{H}_{T,\mathbf{n}}(x) - H(x) \right|
		\leqslant
		\max\left( \left| \widehat{H}^{-}_{T,\mathbf{n}}(x) - H(x) \right|, \left| \widehat{H}^{+}_{T,\mathbf{n}}(x) - H(x) \right| \right).
		\]
		Thus, we conclude that:
		\[
		\left| \widehat{H}_{T,\mathbf{n}}(x) - H(x) \right| = O\left( \left( \frac{k_{T,\mathbf{n}}}{\hat{\mathbf{n}}T} \right)^2 + \sqrt{ \frac{ \log(\hat{\mathbf{n}}T) }{ \hat{\mathbf{n}}T } } \right) \quad \text{a.s.}
		\]
	\end{proof}

	\begin{lemma}\label{lm7}
		Under assumptions \hyperref[H1]{H1}, \hyperref[H2]{H2}, \hyperref[H4]{H4} and \hyperref[H6]{H6}, we have:
		\begin{eqnarray*}
			\sup_{x \in \mathcal{C}}\left|\widehat{H}^{}_{T,\textbf{n}}(x) - H(x)\right| = O\left(\left(\dfrac{k_{T,\textbf{n}}}{\widehat{\textbf{n}}T}\right)^2 +\sqrt{\dfrac{\log(\widehat{\textbf{n}}T)}{\widehat{\textbf{n}}T}}\right) \hspace{1cm}\text{a.s.}
		\end{eqnarray*}
		\label{lem3}
	\end{lemma}
	
	\begin{proof}[Proof of Lemma~\ref{lm7}]
		Since \(D^-_{T,\mathbf{n}}(\beta)\) (resp. \(D^+_{T,\mathbf{n}}(\beta)\)) is non-random, the proof of the almost sure uniform convergence on \(\mathcal{C}\) of \(\widehat{H}^{-}_{T,\mathbf{n}}(x)\) (resp. \(\widehat{H}^{+}_{T,\mathbf{n}}(x)\)) to \(H(x)\) is identical to that given for \(\widetilde{H}_{T,\mathbf{n}}(x)\) in Lemma~\ref{lm2}.
		We choose:
		\[
		D^-_{T,\mathbf{n}}(\beta) = \frac{k_{T,\mathbf{n}}}{\hat{\mathbf{n}}T} \beta^{1/2},
		\quad
		D^+_{T,\mathbf{n}}(\beta) = \frac{k_{T,\mathbf{n}}}{\hat{\mathbf{n}}T} \beta^{-1/2},
		\]
		where \(k_{T,\mathbf{n}}\) is a diverging sequence such that assumption \hyperref[H4]{H4} is satisfied.
		We then obtain:
		\[
		\sup_{x \in \mathcal{C}} \left| \widehat{H}^-_{T,\mathbf{n}}(x) - H(x) \right| = O\left( \left( \frac{k_{T,\mathbf{n}}}{\hat{\mathbf{n}}T} \right)^2 + \sqrt{ \frac{ \log(\hat{\mathbf{n}}T) }{ \hat{\mathbf{n}}T } } \right) \quad \text{a.s.},
		\]
		and similarly:
		\[
		\sup_{x \in \mathcal{C}} \left| \widehat{H}^+_{T,\mathbf{n}}(x) - H(x) \right| = O\left( \left( \frac{k_{T,\mathbf{n}}}{\hat{\mathbf{n}}T} \right)^2 + \sqrt{ \frac{ \log(\hat{\mathbf{n}}T) }{ \hat{\mathbf{n}}T } } \right) \quad \text{a.s.}
		\]
		We have:
		\[
		\left| \widehat{H}_{T,\mathbf{n}}(x) - H(x) \right|  \leqslant  \max\left( \left| \widehat{H}^{-}_{T,\mathbf{n}}(x) - H(x) \right|, \left| \widehat{H}^{+}_{T,\mathbf{n}}(x) - H(x) \right| \right).
		\]
		It follows that:
		\[
		\sup_{x \in \mathcal{C}} \left| \widehat{H}_{T,\mathbf{n}}(x) - H(x) \right| = O\left( \left( \frac{k_{T,\mathbf{n}}}{\hat{\mathbf{n}}T} \right)^2 + \sqrt{ \frac{ \log(\hat{\mathbf{n}}T) }{ \hat{\mathbf{n}}T } } \right) \quad \text{a.s.}
		\]
	\end{proof}
	
	\begin{lemma} \label{lm8}
		Under assumptions \hyperref[H1]{H1}, \hyperref[H2]{H2}, \hyperref[H4]{H4} and \hyperref[H6]{H6}, we have:
		\begin{eqnarray*}
			\left|\widehat{f}^{*}_{T,\textbf{n}}(x) - f^*(x)\right| = O\left(\left(\dfrac{k_{T,\textbf{n}}}{\widehat{\textbf{n}}T} \right)^2+\sqrt{\dfrac{\log \widehat{\textbf{n}}T}{k_{T,\textbf{n}}}}\right) \hspace{1cm} \text{a.s.}
		\end{eqnarray*}
	\end{lemma}
	
	\begin{proof}[Proof of Lemma \ref{lm8}]
		
		Let \(\beta \in ]0,1[\).
		Consider two non-random sequences \(D^-_{T,\mathbf{n}}(\beta)\) and \(D^+_{T,\mathbf{n}}(\beta)\) such that
		\[
		D^-_{T,\mathbf{n}}(\beta) \leqslant  h_{T,\mathbf{n}} \leqslant  D^+_{T,\mathbf{n}}(\beta),
		\]
		where
		\[
		D^-_{T,\mathbf{n}}(\beta) = \left( \dfrac{k_{T,\mathbf{n}}}{\hat{\mathbf{n}}T} \right) \beta^{1/2},
		\quad
		D^+_{T,\mathbf{n}}(\beta) = \left( \dfrac{k_{T,\mathbf{n}}}{\hat{\mathbf{n}}T} \right) \beta^{-1/2}.
		\]
		Note that
		\[
		\dfrac{D^-_{T,\mathbf{n}}(\beta)}{D^+_{T,\mathbf{n}}(\beta)} = \beta.
		\]
		Then, by assumption \hyperref[H2]{H2}, we have
		\begin{align*}
			K\left( \dfrac{x - X_{t,\mathbf{i}}}{h_{T,\mathbf{n}}} \right)
			&= K\left( \dfrac{D^-_{T,\mathbf{n}}(\beta)}{h_{T,\mathbf{n}}} \cdot \dfrac{x - X_{t,\mathbf{i}}}{D^-_{T,\mathbf{n}}(\beta)} \right)
			\geqslant  K\left( \dfrac{x - X_{t,\mathbf{i}}}{D^-_{T,\mathbf{n}}(\beta)} \right),\\
			K\left( \dfrac{x - X_{t,\mathbf{i}}}{D^+_{T,\mathbf{n}}(\beta)} \right)
			&= K\left( \dfrac{h_{T,\mathbf{n}}}{D^+_{T,\mathbf{n}}(\beta)} \cdot \dfrac{x - X_{t,\mathbf{i}}}{h_{T,\mathbf{n}}} \right)
			\geqslant  K\left( \dfrac{x - X_{t,\mathbf{i}}}{h_{T,\mathbf{n}}} \right).
		\end{align*}
		Hence,
		\[
		K\left( \dfrac{x - X_{t,\mathbf{i}}}{D^-_{T,\mathbf{n}}(\beta)} \right)
		\leqslant  K\left( \dfrac{x - X_{t,\mathbf{i}}}{h_{T,\mathbf{n}}} \right)
		\leqslant  K\left( \dfrac{x - X_{t,\mathbf{i}}}{D^+_{T,\mathbf{n}}(\beta)} \right),
		\]
		and also
		\[
		\frac{1}{D^+_{T,\mathbf{n}}(\beta)}  \leqslant  \frac{1}{h_{T,\mathbf{n}}}  \leqslant  \frac{1}{D^-_{T,\mathbf{n}}(\beta)}.
		\]
		It follows that
		\begin{equation}
			\widehat{f}^{-*}_{T,\mathbf{n}}(x)  \leqslant  \widehat{f}^*_{T,\mathbf{n}}(x)  \leqslant  \widehat{f}^{+*}_{T,\mathbf{n}}(x),
			\label{eq13}
		\end{equation}
		where:
		\begin{align*}
			\widehat{f}^{-*}_{T,\mathbf{n}}(x) &= \dfrac{1}{\hat{\mathbf{n}}T D^+_{T,\mathbf{n}}(\beta)} \sum_{\mathbf{i} \in \Lambda_{\mathbf{n}}} \sum_{t=1}^{T} \delta_{t,\mathbf{i}} K\left( \dfrac{x - X_{t,\mathbf{i}}}{D^-_{T,\mathbf{n}}(\beta)} \right),\\
			\widehat{f}^{+*}_{T,\mathbf{n}}(x) &= \dfrac{1}{\hat{\mathbf{n}}T D^-_{T,\mathbf{n}}(\beta)} \sum_{\mathbf{i} \in \Lambda_{\mathbf{n}}} \sum_{t=1}^{T} \delta_{t,\mathbf{i}} K\left( \dfrac{x - X_{t,\mathbf{i}}}{D^+_{T,\mathbf{n}}(\beta)} \right).
		\end{align*}
		Let us now show that \(\widehat{f}^{-*}_{T,\mathbf{n}}\) converges to \(\beta f^*\) for every \(x \in \mathcal{Z}\).
		Indeed,
		\begin{align*}
			\left| \widehat{f}^{-*}_{T,\mathbf{n}}(x) - \beta f^*(x) \right|
			&= \left| \widehat{f}^{-*}_{T,\mathbf{n}}(x) - \mathbb{E}[\widehat{f}^{-*}_{T,\mathbf{n}}(x)] + \mathbb{E}[\widehat{f}^{-*}_{T,\mathbf{n}}(x)] - \beta f^*(x) \right|\\
			& \leqslant  \left| \widehat{f}^{-*}_{T,\mathbf{n}}(x) - \mathbb{E}[\widehat{f}^{-*}_{T,\mathbf{n}}(x)] \right| + \left| \mathbb{E}[\widehat{f}^{-*}_{T,\mathbf{n}}(x)] - \beta f^*(x) \right|.
		\end{align*}
		Following the same steps as in the proofs of Lemmas~\ref{lm2} and~\ref{lm4}, we obtain
		\[
		\left| \widehat{f}^{-*}_{T,\mathbf{n}}(x) - \beta f^*(x) \right| = O\left( \left( \frac{k_{T,\mathbf{n}}}{\hat{\mathbf{n}}T} \right)^2 + \sqrt{ \frac{\log(\hat{\mathbf{n}}T)}{k_{T,\mathbf{n}}} } \right) \quad \text{a.s.}
		\]
		Similarly, we have:
		\[
		\left| \widehat{f}^{+*}_{T,\mathbf{n}}(x) - \frac{1}{\beta} f^*(x) \right| = O\left( \left( \frac{k_{T,\mathbf{n}}}{\hat{\mathbf{n}}T} \right)^2 + \sqrt{ \frac{\log(\hat{\mathbf{n}}T)}{k_{T,\mathbf{n}}} } \right) \quad \text{a.s.}
		\]
		Defining
		\[
		u_{T,\mathbf{n}} = \left( \frac{k_{T,\mathbf{n}}}{\hat{\mathbf{n}}T} \right)^2 + \sqrt{ \frac{ \log(\hat{\mathbf{n}}T) }{ k_{T,\mathbf{n}} } },
		\quad
		\beta = \beta_{T,\mathbf{n},\epsilon} = 1 - \frac{ \epsilon u_{T,\mathbf{n}} }{ 3f^* },
		\quad \forall \epsilon > 0,
		\]
		we deduce that
		\[
		\left| \widehat{f}^{-*}_{T,\mathbf{n}}(x) - f^*(x) \right| = O(u_{T,\mathbf{n}}) \quad \text{a.s.}
		\]
		Indeed,
		\begin{align*}
			\left| \widehat{f}^{-*}_{T,\mathbf{n}}(x) - f^*(x) \right|
			&= \left| \widehat{f}^{-*}_{T,\mathbf{n}}(x) - \beta f^*(x) + \beta f^*(x) - f^*(x) \right|\\
			& \leqslant  \left| \widehat{f}^{-*}_{T,\mathbf{n}}(x) - \beta f^*(x) \right| + \left| - \frac{ \epsilon u_{T,\mathbf{n}} }{3} \right|.
		\end{align*}
		Similarly, for \(\widehat{f}^{+*}_{T,\mathbf{n}}\), we have
		\begin{align*}
			\left| \widehat{f}^{+*}_{T,\mathbf{n}}(x) - f^*(x) \right|
			&= \left| \widehat{f}^{+*}_{T,\mathbf{n}}(x) - \frac{1}{\beta} f^*(x) + \frac{1}{\beta} f^*(x) - f^*(x) \right|\\
			& \leqslant  \left| \widehat{f}^{+*}_{T,\mathbf{n}}(x) - \frac{1}{\beta} f^*(x) \right| + \left| \frac{ \epsilon u_{T,\mathbf{n}} }{3} \right|.
		\end{align*}
		A Taylor expansion at 0 makes it possible to verify that \(\beta_{T,\mathbf{n},\epsilon} - 1\) and \(\frac{1}{\beta_{T,\mathbf{n},\epsilon}} - 1\) converge to 0 at the same rate.
		Finally, from
		\[
		\widehat{f}^{-*}_{T,\mathbf{n}}(x)  \leqslant  \widehat{f}^*_{T,\mathbf{n}}(x)  \leqslant  \widehat{f}^{+*}_{T,\mathbf{n}}(x),
		\]
		which implies
		\[
		\left| \widehat{f}^*_{T,\mathbf{n}}(x) - f^*(x) \right|  \leqslant  \max\left( \left| \widehat{f}^{-*}_{T,\mathbf{n}}(x) - f^*(x) \right|, \left| \widehat{f}^{+*}_{T,\mathbf{n}}(x) - f^*(x) \right| \right),
		\]
		we obtain
		\[
		\left| \widehat{f}^*_{T,\mathbf{n}}(x) - f^*(x) \right| = O\left( \left( \frac{k_{T,\mathbf{n}}}{\hat{\mathbf{n}}T} \right)^2 + \sqrt{ \frac{ \log(\hat{\mathbf{n}}T) }{ k_{T,\mathbf{n}} } } \right) \quad \text{a.s.}
		\]
		
	\end{proof}
	
	\begin{lemma}\label{lm9}
		Under assumptions \hyperref[H1]{H1}, \hyperref[H2]{H2}, \hyperref[H4]{H4} and \hyperref[H6]{H6}, we have:
		\begin{eqnarray*}
			\sup_{x \in \mathcal{C}}\left|\widehat{f}^{*}_{T,\textbf{n}}(x) - f^*(x)\right| = O\left(\left(\dfrac{k_{T,\textbf{n}}}{\widehat{\textbf{n}}T} \right)^2+\sqrt{\dfrac{\log \widehat{\textbf{n}}T}{k_{T,\textbf{n}}}}\right) \hspace{1cm} \text{a.s.}
		\end{eqnarray*}
	\end{lemma}
	\begin{proof}[Proof of Lemma \ref{lm9}]
		By using the same reasoning as in the proof of Lemma~\ref{lm8}, and following the approach adopted in the proofs of Lemmas~\ref{lm3} and~\ref{lm5}, we obtain the announced result.
	\end{proof}

\subsection{Proof of Theorem~\ref{th1}}

\begin{proof}[Proof of Theorem \ref{th1}]
		
		For every \(x \in \mathcal{Z}\), consider the following decomposition:
		\[
		\left|\widetilde{\lambda}_{T,\mathbf{n}}(x) - \lambda(x)\right|
		\;\leqslant \;
		\frac{1}{1 - \widetilde{H}_{T,\mathbf{n}}(x)}
		\left[
		\left|\widetilde{f}^*_{T,\mathbf{n}}(x) - f^*(x)\right|
		+
		\left|\frac{f^*(x)}{1 - H(x)}\right|
		\left|\widetilde{H}_{T,\mathbf{n}}(x) - H(x)\right|
		\right].
		\]
		There exists a constant \(\mu > 0\) such that, for every \(x \in \mathcal{Z}\),
		\[
		H(x) \leqslant  1 - \mu,
		\quad \text{which implies} \quad
		1 - H(x) \geqslant  \mu.
		\]
		Moreover, we know that \(\widetilde{H}_{T,\mathbf{n}}(x)\) converges almost surely to \(H(x)\); in other words, there exists an event \(\Omega_0\) with probability 1 such that, for every \(\omega \in \Omega_0\) and every \(x \in \mathcal{Z}\),
		\[
		\left| \widetilde{H}_{T,\mathbf{n}}(\omega,x) - H(x) \right| \xrightarrow[\;(T,\mathbf{n}) \to +\infty\;]{} 0.
		\]
		It follows that, for every \(\varepsilon > 0\) with \(\varepsilon < \mu\), there exists an integer \(m(\omega)\) (independent of \(x\)) such that for all \(T, \hat{\mathbf{n}} \geqslant m(\omega)\),
		\[
		\widetilde{H}_{T,\mathbf{n}}(\omega,x)
		\;\leqslant \;
		H(x) + \varepsilon
		\;\leqslant \;
		(1 - \mu) + \varepsilon
		= 1 - (\mu - \varepsilon),
		\]
		which implies:
		\[
		1 - \widetilde{H}_{T,\mathbf{n}}(x)
		\geqslant
		\mu - \varepsilon > 0
		\quad \text{(almost surely)}.
		\]
		Consequently, for sufficiently large \(T\) and \(\hat{\mathbf{n}}\),
		\[
		\frac{1}{1 - \widetilde{H}_{T,\mathbf{n}}(x)}
		\;\underset{\text{a.s.}}{\leqslant }\;
		\frac{1}{\mu - \varepsilon}.
		\]
		Furthermore, under assumption \hyperref[H3]{H3}, we have:
		\[
		\left| \frac{f^*(x)}{1 - H(x)} \right|
		\;\leqslant \;
		\frac{\upsilon}{\mu},
		\]
		where \(\upsilon > 0\) is such that \(f^*(x) \leqslant  \upsilon\) for every \(x \in \mathcal{Z}\).
		Combining the previous inequalities, we deduce that there exists an event with probability 1 (the same \(\Omega_0\)) such that, for sufficiently large \(T\) and \(\hat{\mathbf{n}}\), and for every \(x \in \mathcal{Z}\),
		\begin{equation}
			\left|\widetilde{\lambda}_{T,\mathbf{n}}(x) - \lambda(x)\right|
			\;\underset{\text{a.s.}}{\leqslant }\;
			\frac{1}{\mu - \varepsilon}
			\left[
			\left|\widetilde{f}^*_{T,\mathbf{n}}(x) - f^*(x)\right|
			+
			\frac{\upsilon}{\mu}
			\left| \widetilde{H}_{T,\mathbf{n}}(x) - H(x) \right|
			\right].
			\label{eq8}
		\end{equation}
		Finally, using the results of Lemmas~\ref{lm2} and~\ref{lm4}, we obtain:
		\[
		\left| \widetilde{\lambda}_{T,\mathbf{n}}(x) - \lambda(x) \right|
		=
		O\left(
		\widetilde{h}_{T,\mathbf{n}}^2
		+
		\sqrt{ \frac{ \log(\hat{\mathbf{n}}T) }{ \hat{\mathbf{n}}T \widetilde{h}_{T,\mathbf{n}} } }
		\right)
		\quad \text{a.s.}
		\]
	\end{proof}

\subsection{Proof of Theorem~\ref{th2}}

\begin{proof}[Proof of Theorem \ref{th2}]
		
		Let \(x\in\mathcal C\). By \hyperref[H7]{H7}, there exists \(\mu>0\) such that
		$
		\inf_{x\in\mathcal C}\{1-H(x)\}\geqslant \mu.
		$
		This uniform bound on \(\mathcal C\) allows us to take the supremum over \(x\) in \eqref{eq8}. Thus, on the event \(\Omega_0\) of probability \(1\), we obtain:
		\[
		\sup_{x \in \mathcal{C}} \left| \widetilde{\lambda}_{T,\mathbf{n}}(x) - \lambda(x) \right|
		\;\underset{\text{a.s.}}{\leqslant }\;
		\frac{1}{\mu - \varepsilon}
		\left(
		\sup_{x \in \mathcal{C}} \left| \widetilde{f}^*_{T,\mathbf{n}}(x) - f^*(x) \right|
		+
		\frac{\upsilon}{\mu}
		\sup_{x \in \mathcal{C}} \left| \widetilde{H}_{T,\mathbf{n}}(x) - H(x) \right|
		\right).
		\]
		Applying next the results of Lemmas~\ref{l3} and~\ref{l5}, which provide respectively the uniform convergence rates of \(\widetilde{f}^*_{T,\mathbf{n}}\) to \(f^*\) and of \(\widetilde{H}_{T,\mathbf{n}}\) to \(H\), it follows that:
		\[
		\sup_{x \in \mathcal{C}} \left| \widetilde{\lambda}_{T,\mathbf{n}}(x) - \lambda(x) \right|
		=
		O\left(
		\widetilde{h}_{T,\mathbf{n}}^2
		+
		\sqrt{ \frac{ \log(\hat{\mathbf{n}}T) }{ \hat{\mathbf{n}}T \widetilde{h}_{T,\mathbf{n}} } }
		\right)
		\quad \text{a.s.}
		\]
		
	\end{proof}

	\section*{Declarations}
	
	\noindent
	\textbf{Funding} \\[2pt]
	No specific funding was received for this research.
	
	\vspace{0.5em}
	\noindent
	\textbf{Competing interests} \\[2pt]
	The authors declare no competing interests.
	
	\vspace{0.5em}
	\noindent
	\textbf{Author contributions} \\[2pt]
	All authors contributed equally to all aspects of this work.

\bibliographystyle{plain}
\bibliography{references}

\end{document}